\documentclass[a4paper,12pt]{amsart}
\usepackage{color}
\usepackage{amsmath,amssymb,amscd,amsfonts,amsthm,stmaryrd}
\usepackage[mathscr]{eucal}
\usepackage{xcolor}
\usepackage{graphicx}
\usepackage{tikz-cd}
\usepackage{adjustbox}

\usepackage{nicefrac}

\usepackage{color}

\numberwithin{equation}{section}

\def\diam{\operatorname{diam}}
\def\rank{\operatorname{rank}}
\def\pol{\operatorname{pol}}
\def\rad{\operatorname{rad}}

\def\Ant{\operatorname{Ant}}

\def\acts{\curvearrowright}
\def\D{\partial}

\def\R{\mathbb R}
\def\Z{\mathbb Z}

\def\N{\mathbb N}

\def\al{\alpha}

\def\ka{\kappa}

\def\eps{\epsilon}
\def\ga{\gamma}
\def\Ga{\Gamma}

\def\la{\lambda}

\def\om{\omega}

\def\si{\sigma}
\def\Si{\Sigma}

\def\tits{\partial_{T}}
\def\geo{\partial_{\infty}}

\def\wlim{\mathop{\hbox{$\om$-lim}}}

\def\acts{\curvearrowright}

\def\<{\langle}
\def\>{\rangle}

\theoremstyle{plain}
\newtheorem{thm}{Theorem}[section]
\newtheorem{lem}[thm]{Lemma}
\newtheorem{prop}[thm]{Proposition}
\newtheorem{cor}[thm]{Corollary}
\newtheorem{slem}[thm]{Sublemma}

\newtheorem{introthm}{Theorem}

\newtheorem{introcor}[introthm]{Corollary}

\newtheorem*{conj}{Conjecture}

\newtheoremstyle{named}{}{}{\itshape}{}{\bfseries}{.}{.5em}{\thmnote{#3} #1}
\theoremstyle{named}
\newtheorem*{namedthm}{Theorem}

\theoremstyle{definition}
\newtheorem{dfn}[thm]{Definition}

\theoremstyle{remark}
\newtheorem{rem}[thm]{Remark}

\newcommand{\bcl}{\begin{claim}}
\newcommand{\ecl}{\end{claim}}
\newcommand{\bcor}{\begin{cor}}
\newcommand{\ecor}{\end{cor}}
\newcommand{\bdfn}{\begin{dfn}}
\newcommand{\edfn}{\end{dfn}}
\newcommand{\ben}{\begin{enumerate}}
\newcommand{\bit}{\begin{itemize}}
\newcommand{\blem}{\begin{lem}}
\newcommand{\bslem}{\begin{slem}}
\newcommand{\bprop}{\begin{prop}}
\newcommand{\bthm}{\begin{thm}}
\newcommand{\een}{\end{enumerate}}
\newcommand{\eit}{\end{itemize}}
\newcommand{\elem}{\end{lem}}
\newcommand{\eslem}{\end{slem}}
\newcommand{\eprop}{\end{prop}}
\newcommand{\ethm}{\end{thm}}

\usepackage{hyperref}



\begin{document}

\title[Higher rank II]{CAT(0) spaces of higher rank II}
\author{Stephan Stadler}

\newcommand{\Addresses}{{\bigskip\footnotesize
\noindent Stephan Stadler,
\par\nopagebreak\noindent\textsc{Max Planck Institute for Mathematics, Vivatsgasse 7, 53111 Bonn, Germany}
\par\nopagebreak
\noindent\textit{Email}: \texttt{stadler@mpim-bonn.mpg.de}

}}


\begin{abstract}
This belongs to a series of papers 
motivated by Ballmann's  Higher Rank Rigidity Conjecture. We prove the following.
Let $X$ be a CAT(0) space  with a geometric group action $\Ga\acts X$.  
Suppose that every geodesic in $X$ lies in an $n$-flat, $n\geq 2$. 
If $X$ contains a periodic $n$-flat which does not bound a flat $(n+1)$-half-space, then $X$ is a Riemannian symmetric space, a Euclidean building 
or  non-trivially splits as a metric product.
This generalizes the Higher Rank Rigidity Theorem for Hadamard manifolds with geometric group actions.
\end{abstract}

\maketitle

\tableofcontents

\section{Introduction}

\subsection{Main results}
A central role in the geometry of nonpositively curved spaces is played by {\em flats} -- isometric embeddings of Euclidean spaces.
One expects that in a CAT(0) space with enough extremal curvature zero, flats are uniformly distributed throughout the space and special geometry appears --
at least if the isometry group is sufficiently large. This is made precise by Ballmann's conjecture which is the main motivation for the present paper.

\begin{conj}[Higher Rank Rigidity]\label{conj_rr}
Let $X$ be a locally compact CAT(0) space with a geometric group action $\Ga\acts X$.
If every geodesic in $X$ lies in an $n$-flat, $n\geq 2$, then $X$ is a Riemannian symmetric space or a Euclidean
building or non-trivially splits as a metric product.
\end{conj}

Recall that an isometric action is called {\em geometric} if it is properly discontinuous and cocompact. 
For Hadamard manifolds Higher Rank Rigidity is a celebrated theorem by Ballmann \cite{B_higher} and independently Burns-Spatzier \cite{BS_higher}, see also
Eberlein-Heber \cite{EH_diff}.
%
%

In \cite{St_rrIII} we have confirmed Higher Rank Rigidity as well as the Closing Lemma for spaces without 3-flats.
More precisely, we obtained

\bthm[{\cite[Theorems~A and~B]{St_rrIII}}]\label{thm_rr}
Let $X$ be a locally compact  CAT(0) space without 3-flats. 
Suppose $X$ admits a geometric group action $\Ga\acts X$. If $X$ contains a complete geodesic which does not bound a flat half-plane, then it also contains a
 $\Ga$-periodic geodesic which does not bound a flat half-plane. On the other hand, if every complete geodesic in $X$ bounds a flat half-plane,
then $X$ is a Riemannian symmetric space, a Euclidean building or non-trivially splits as a metric product.
\ethm

Moreover, we have proved the following general result.

\bthm[{\cite[Theorem~A]{St_rr}}]\label{thm_rigidity}
Let $X$ be a locally compact  CAT(0) space whose Tits boundary has dimension $n-1\geq 1$. 
Suppose that every geodesic in $X$ lies in an $n$-flat. 
If $X$ contains a periodic $n$-flat, then $X$ is a Riemannian symmetric space or a Euclidean building or  $X$ non-trivially splits as a metric product.
\ethm

A flat is periodic, if its stabilizer in the isometry group contains a subgroup acting geometrically on the flat.
Recall that if $X$ is a rank $n$ symmetric space or  Euclidean building with a geometric group action, then periodic $n$-flats are dense in the space of $n$-flats in $X$ \cite[Theorem~8.9]{BaBr_orbi}, \cite[Lemma~8.3]{Mos_rigidity}, \cite[Theorem~2.8]{PR_cartan}. Note that Theorem~\ref{thm_rigidity} does not make any symmetry assumptions besides the periodic flat.

We say that a CAT(0) space has {\em rank at least $n$},  $\rank(X)\geq n$, if every geodesic in $X$ lies in an $n$-flat\footnote{
See \cite{OR_int} for a definition of rank for CAT(0) spaces whose isometry groups satisfy the duality condition.}.
In the present paper we consider CAT(0) spaces $X$ of {\em higher rank}, i.e. rank at least $2$. We
 will work under the assumption of a geometric group action but remove any condition on the dimension of the Tits boundary.
Note  that the conditions $\rank(X)=n\geq 2$ and $\dim(\tits X)=n-1$
in combination ensure that the Tits boundary is geodesically complete (\cite[Lemma~2.1]{BL_building}).
This allowed us in \cite{St_rr} to rely on a deep result of Lytchak:

\bthm[{\cite[Main Theorem]{Ly_rigidity}}]\label{thm_lytchak}
Let $Z$ be a finite-dimensional geodesically complete
CAT(1) space. Suppose $Z$ has a proper closed subset $A$ which is {\em symmetric}, i.e. contains with each
point all of its antipodes. Then $Z$ is
a spherical join or a spherical building.
\ethm

If the Tits boundary of a locally compact geodesically complete CAT(0) space is a spherical building, then rigidity follows from Leeb's theorem \cite[Main Theorem]{Leeb}.
However, Tits boundaries of CAT(0) spaces tend to be not geodesically complete.
One can find examples even among $4$-dimensional real-analytic Hadamard manifolds which admit geometric group actions \cite{HuSch_tits}.
If $V$ is such a manifold, then $V$ has necessarily rank $1$. 
However, the product with any other CAT(0) space $X$ has higher rank, yet $V\times X$ cannot be subject to Theorem~\ref{thm_rigidity}
as it has to contain maximal flats of different dimensions. Vice versa, if Higher Rank Rigidity holds, then a higher rank space can only contain
maximal flats of different dimensions if it splits non-trivially.

Recall that a flat is called {\em Morse}, if it does not bound a flat half-space.

\begin{namedthm}[Main]\label{thm_mainB}
Let $X$ be a locally compact CAT(0) space of rank at least $n\geq 2$ with a geometric group action $\Ga\acts X$.  
If $X$ contains a periodic Morse $n$-flat, then $X$ is a Riemannian symmetric space or a Euclidean building or  $X$
non-trivially splits as a metric product.
\end{namedthm}

\begin{introcor}\label{cor_mainA}
Let $M$ be a Hadamard manifold with a geometric group action $\Ga\acts M$.
If $M$ has higher rank, then $M$ is a Riemannian symmetric space or non-trivially splits as a metric product.
\end{introcor}



\subsection{Strategy}

The proof of our \hyperref[thm_mainB]{Main Theorem} relies on Leeb's rigidity theorem.

\bthm[{\cite[Main Theorem]{Leeb}}]\label{thm_leeb}
Let $X$ be a locally compact, geodesically complete CAT(0) space. If $\tits X$ is a connected thick irreducible spherical building,
then $X$ is either a symmetric space or a Euclidean building.
\ethm

As we have explained above, in our setting the Tits boundary might not be geodesically complete.
Thus we cannot use Lytchak's rigidity theorem to show that the Tits boundary is a spherical join or a spherical building.
Our starting point is the existence of regular points at infinity: 

\bthm[{\cite[Theorem~B]{St_rr}}]\label{thm_regular} 
Let $X$ be a locally compact CAT(0) space.  Suppose that $X$ contains a periodic $n$-flat $F$. Then either $F$ bounds a flat $(n+1)$-half-space, or 
$\tits F\subset\tits X$ contains a dense subset which is open in $\tits X$. 
\ethm

Using dynamical ideas based on work of Papasoglu-Swenson \cite{PaSw_bound} and Guralnik-Swenson \cite{GuSw_trans}, 
we then show the existence of a non-empty proper closed symmetric subset in $\geo X$.
Note that this is not the same setup than in Theorem~\ref{thm_lytchak}.
It is essential, that we find such sets which are closed in the cone topology and not only in the Tits topology.
The proof then relies on the basic observation that 
minimal closed symmetric subsets constitute an equidistant decomposition of the Tits boundary, and therefore, 
a well-defined metric quotient -- a submetry $\delta:\tits X\to\Delta$.
We then study this submetry carefully.
First control on $\delta$ comes from the fact that
each logarithmic map from the Tits boundary to the space of directions at some point defines
a submetry between the corresponding quotients. 
Using logarithmic maps for varying sets of points,  we can
show that the Tits boundary is geodesically complete in horizontal directions.
A key role is then played by the periodic Morse flat $\hat F$.
It turns out that $\delta$ restricts to a submetry $\hat\delta:\tits \hat F\to\Delta$. 
We show that if $\tits X$ is irreducible, then $\hat\delta$
is discrete and $\Delta$ is the quotient of $\tits \hat F$ by a finite group action. We conclude the proof by showing that $\tits \hat F$
has to be top-dimensional if $\tits X$ is irreducible. In this case Theorem~\ref{thm_rigidity} from \cite{St_rr} applies.
The last step relies on the fact that submetries of round spheres to spaces of constant curvature 1 are very well understood, in fact classified.

\medskip

\subsection{Organization}
In Section~\ref{sec_pre} we recall necessary basics from metric geometry and agree on notation.
Apart from several standard facts, we also review the structure of submetries from round spheres and prove several results  needed later on. 
We also recall the notion of a Morse flat and its basic properties.
In Section~\ref{sec_sym} we introduce symmetric and visually symmetric subsets,  preparing for the structure of a submetry.
We also proof a first basic distance estimate induced by proper closed symmetric subsets (Lemma~\ref{lem_bound}) which is ultimately responsible 
for rigidity.
Further, we discuss almost near points. Segments between almost near points always extend to geodesics of length $\pi$. 
This will be used to show geodesic completeness in horizontal directions of submetries. We also show that almost near points lead to many spherical triangles.
In Section~\ref{sec_sub} we show that closed minimal symmetric sets provide a submetry of the Tits boundary (Proposition~\ref{prop_submetry}).
We then discuss the branching behavior of horizontal geodesics and prove splitting results based on particular splittings of the Tits boundary of a
periodic Morse flat.
We begin Section~\ref{sec_final} by showing that a periodic Morse flat in a higher rank space leads to a non-trivial proper symmetric subset at infinity
and therefore to a submetry of the Tits boundary. We then study the Tits boundary of a periodic Morse flat within the Tits boundary of our space and show
that it carries a structure reminiscent of a Coxeter complex. We then consider induced submetries and show in the irreducible case that the dimension of
the base of our submetry is equal to the dimension of the Tits boundary of the periodic Morse flat.
We prove the \hyperref[thm_mainB]{Main Theorem} by showing that the dimension of the base is equal to the dimension of the Tits boundary unless the space is reducible.
At last, we conclude Corollary~\ref{cor_mainA}.

\subsection{Acknowledgments}
It's my pleasure to thank Alexander Lytchak for all his support.
This article owes a lot to his paper \cite{Ly_rigidity}. In addition I want to thank him for
many helpful discussions, in particular on submetries. I also want to thank Bruce Kleiner and Bernhard Leeb
for inspiring discussions. 
I was supported by DFG grant SPP 2026.

\section{Preliminaries}\label{sec_pre}

References for this section are \cite{AKP, Ballmann, BH, KleinerLeeb, Ly_rigidity}.

\subsection{Metric spaces}
Euclidean $n$-space with its flat metric will be denoted by $\R^n$ and the unit sphere by $S^{n-1}\subset\R^n$.
We denote  the distance between two points $x$ and $y$ in a metric space $X$ by $|x,y|$.
If $A\subset X$ denotes a subset, then $|x,A|$ refers to the greatest lower bound for distances from points in $A$ to $x$.
Similarly, for subsets $A, A'\subset X$ we denote by $|A,A'|$ the greatest lower bound for distances from points in $A$ to points in $A'$.  
For $x\in X$ and $r>0$, we denote by $B_r(x)$ and $\bar B_r(x)$ the open and closed $r$-ball around $x$, respectively.
A \emph{geodesic}
is an isometric embedding of an interval. It is called a {\em geodesic segment}, if it is compact.
The {\em endpoints} of a geodesic segment $c$ are denoted by $\D c$.
A geodesic segment $c$ {\em branches} at an endpoint $y\in \D c=\{x,y\}$, if there are geodesics $c^\pm$ starting in $x$ which  strictly contain $c$
and such that $c^-\cap c^+=c$.

 A \emph{triangle} is a union of three geodesics connecting three points.
If $x,y,z$ are point with unique geodesics between them, then we denote the corresponding triangle by $\triangle(x,y,z)$.
$X$  is a \emph{geodesic metric space} if
any pair of   points of $X$
is connected by a geodesic.
It is \emph{geodesically complete} if every geodesic segment is contained in a complete local geodesic.

For a sequence of metric spaces $(X_k)$ with uniformly bounded diameter, we denote by $X_\om:=\wlim X_k$ its ultralimit
with respect to a non-principal ultrafilter $\om$ on the natural numbers.
Recall that if a sequence of surjective 1-Lipschitz maps $f_k:X_k\to Y_k$ between metric spaces $X_k, Y_k$
with uniformly bounded diameters induce a surjective 1-Lipschitz maps $f_\om:X_\om\to Y_\om$.

\subsection{Submetries}

In this section we recall the notion of a submetry. Originally the concept was introduced by Berestovskii in \cite{Ber_sub}
as a metric analog of a Riemannian submersion. If $B$ is a spherical building, then the natural $1$-Lipschitz map $\delta: B\to\Delta_{mod}$ 
which folds the building onto is model chamber is a submetry. Vice versa,  CAT(1) spaces
which admit certain kinds of submetries are spherical buildings \cite{Ly_rigidity}.
This explains the relevance of submetries for us. 
\medskip

Two subsets $X^\pm$ in a metric space $X$ are called {\em equidistant}, if for every point $x^\pm\in X^\pm$
there exists a point $y^\mp\in X^\mp$ such that $|x^\pm,y^\mp|=|X^\pm,X^\mp|$ holds.

A decomposition of a metric space $X$ into closed pairwise equidistant subsets $\{X_y\subset X|\ y\in Y\}$ induces a natural metric on the quotient space $Y$:
\[|y,y'|:=|X_y,X_{y'}|.\]
The natural projection $\delta:X\to Y$ is a submetry in the sense of \cite{Ber_sub}.
\bdfn
A {\em submetry} $\delta:X\to Y$ between metric spaces is a map which sends for every point $x$ in $X$ and each radius $r\geq 0$
the closed $r$-ball around $x$ surjectively onto the closed $r$-ball around $\delta(x)$.
\edfn
Submetries on a metric space $X$ are in one-to-one correspondence to closed equidistant decompositions of $X$.
Each submetry is a 1-Lipschitz map. The composition of two submetries is again a submetry.
Vice versa, if two submetries $f:X\to Y$ and $h:X\to Z$ have the property that every $f$-fiber is contained in an $h$-fiber,
then the natural map  $g:Y\to Z$ with $g\circ f=h$ is a submetry.

If $\delta:X\to Y$ is a submetry, then
two points $x,x'$ in $X$ are called $\delta$-near, if $|x,x'|=|\delta(x),\delta(x')|$ holds.
In this case, $\delta$ maps any geodesic between $x$ and $x'$ isometrically onto a geodesic between $\delta(x)$
and $\delta(x')$. Such geodesics in $X$ are called {\em horizontal}. More generally, a geodesic in $X$ is called {\em piecewise horizontal},
if it is the union of finitely many horizontal geodesics. 

We collect a couple of results on submetries of round spheres which will be used later on.
Most of these results hold in greater generality but we restrict to the cases of interest for us.

Let $\delta:S^n\to\Delta$ be a submetry. Then $\Delta$ is an Alexandrov space dimension $k\leq n$ and curvature at least 1 \cite[Proposition~3.1]{KaLy_sub},
\cite{BGP}.
For any rectifiable curve $\ga$ in $\Delta$ starting in a point $\delta(x)$ we find a rectifiable curve $\tilde\ga$ to $S^n$ starting in $x$ such that 
$\delta\circ\tilde\ga=\ga$ and $\delta$ preserves the length of $\tilde \ga$ \cite[Lemma~2.8]{KaLy_sub}. Such a curve is called a
{\em horizontal lift} of $\ga$.
The submetry $\delta$ is differentiable  at every point $x\in S^n$ in the following sense. 
There exists a map $d\delta(x):T_x S^n\to T_{\delta(x)}\Delta$, such that for every
sequence $t_i\to 0$, the submetry $\delta$ seen as a map between rescaled spaces
$\delta:(\frac{1}{t_i}S^n,x)\to(\frac{1}{t_i}\Delta,\delta(x))$ converge to the map $d\delta(x)$.
The differential $d\delta(x):T_x S^n\to T_{\delta(x)}\Delta$
is an {\em infinitesimal submetry}, a submetry which commutes with the natural dilations on the Euclidean cones 
$T_x S^n$ and $T_{\delta(x)}\Delta$ \cite[Proposition~3.3]{KaLy_sub}.

The base $\Delta$ allows for a natural disjoint decomposition into strata $\Delta^l\subset\Delta$, $0\leq l\leq k$.
The $l$-dimensional stratum $\Delta^l$ is
the set of points $y\in\Delta$ whose tangent space $T_y\Delta$ splits off  $\R^l$ as a direct
factor but not $\R^{l+1}$.

\bthm[{\cite[Theorem~1.6]{KaLy_sub}}]
Let $\delta:S^n\to\Delta$ be a submetry with $\dim(\Delta)=k$. For any $0\leq l\leq k$, the stratum $\Delta^l$ is
an $l$-dimensional topological manifold which is locally closed and locally
convex in $\Delta$. The maximal stratum $\Delta^k$ is open and globally convex.
\ethm

If a geodesic $\ga$ in $\Delta$ starts in a point of $\Delta^l_+:=\bigcup_{i=l}^k \Delta^i$, then it stays in this set \cite[Lemma~10.1]{KaLy_sub}.
Thus the codimension 1 stratum $\Delta^{k-1}$ consists precisely of the points $y\in\Delta$ where the tangent space $T_y\Delta$
is isometric to a $k$-dimensional Euclidean half-space. The top-dimensional stratum is sometimes referred to as the set of {\em regular points},
$\Delta_{reg}:=\Delta^k$.

\bthm[{\cite[Theorem~10.5]{KaLy_sub}}]\label{thm_submersion}
Let $\delta:S^n\to\Delta$ be a submetry. Then $\delta:\delta^{-1}(\Delta_{reg})\to\Delta_{reg}$
is a $\mathcal{C}^{1,1}$ Riemannian submersion. 
\ethm

We call a submetry $\delta:S^n\to\Delta$ {\em transnormal}, if any extension of a $\delta$-horizontal geodesic in $S^n$ is a piecewise horizontal
geodesic.
This is equivalent to saying that an antipode of a direction tangent to a $\delta$-horizontal geodesic is again tangent to a $\delta$-horizontal geodesic.
Since the horizontal space $H_x\subset\Si_x S^n$ is convex for any submetry \cite[Proposition~12.5]{KaLy_sub}, it follows that for transnormal submetries
horizontal spaces are round spheres. We even have:

\bprop[{\cite[Proposition~12.5]{KaLy_sub}}]\label{prop_C11}
Let $\delta:S^n\to\Delta$ be a transnormal submetry. Then every fiber of $\delta$ is a $\mathcal{C}^{1,1}$ submanifold in $S^n$.
\eprop


\blem\label{lem_polytope}
Let $\delta:S^n\to\Delta$ be a submetry such that every geodesic in $S^n$ is piecewise $\delta$-horizontal. Then
$\Delta$ is a spherical orbifold and $\delta$ is a covering of Riemannian orbifolds.
\elem

\proof
By \cite[Proposition~9.1]{Ly_diss}, the regular fibers of $\delta$ are discrete and $\dim(\Delta)=n$.
The claim follows from \cite[Theorem~1.2]{La_orbi}.
\qed

%

\blem\label{lem_sphhemi}
Let $\delta:S^m\to S^n_+$ be a submetry where $S^n_+\subset S^n$ is a hemisphere.
Then $\delta$ splits as $(\delta_1,\delta_2):S^{n-1}\circ S^{m-n}\to S^{n-1}\circ\{v\}$
where $\delta_1$ is an isometry and $\delta_2$ is constant.
\elem

\proof
Write $S^{n}_+=\{v_1^+,v_1^-\}\circ\ldots\circ\{v_n^+,v_n^-\}\circ\{v\}$ and observe that the inverse image under $\delta$ of
a pair of antipodes $\{v_j^+,v_j^-\}$ has to be a pair of antipodes. The claim follows.
\qed

\blem\label{lem_fold}
Let $\delta:S^m\to\Delta$ be a transnormal submetry.
Let $x\in S^m$ be a point with $\delta(x)=y\in\Delta^{n-1}$. Then the horizontal space $H_x\subset\Si_x S^m$
admits a splitting $H_x\cong K\circ K^\perp$ such that the differential $d\delta(x)$
maps $K$ isometrically onto $\Si_y\Delta^{n-1}\cong S^{n-2}$ and sends $K^\perp$ to the point $w\in\Si_y\Delta$ where $\Si_y\Delta\cong\Si_y \Delta^{n-1}\circ\{w\}$.
In case $n=1$, $K$ is empty.
\elem

\proof
We may assume $n\geq 2$, the case $n=1$ will then be clear.
Since $\delta$ is transnormal, the horizontal space $H_x\subset\Si_x S^m$ is a round sphere.
By \cite[Proposition~5.5]{KaLy_sub}, the differential $d\delta(x)$ restricts to a submetry $d\delta^H(x):H_x\to\Si_y\Delta$.
Recall that $\Si_y\Delta\cong\Si_y\Delta^{n-1}\circ\{w\}$ and $\Si_y\Delta^{n-1}\cong S^{n-2}$. 
The claim follows from Lemma~\ref{lem_sphhemi}.
\qed

\blem\label{lem_focal}
Let $\delta:S^m\to\Delta$ be a transnormal submetry.
Suppose that $\Delta$ is a spherical orbifold of dimension $n<m$.
Then $\Delta^{n-1}$ is non-empty and there exists a point $x\in S^m$ with $\delta(x)\in\Delta^{n-1}$ such that 
the horizontal space $H_x\subset\Si_x S^m$ is a round sphere of dimension at least $n$. 
\elem

This follows rather quickly form the general theory, as has been explained to us by Alexander Lytchak.
\proof By \cite{MR_lap} it is enough to consider the case where $\delta$ has connected smooth fibers.
Then, the restriction $\delta:\delta^{-1}(\Delta_{reg})\to\Delta_{reg}$ is a smooth Riemannian submersion.
Since the curvature on $\Delta_{reg}$ is equal to the curvature of $S^n$, the horizontal distribution is integrable.
As a consequence, all regular fibers are isoparametric submanifolds of $S^n$ \cite[Definition~2.1]{Thor}.
Thus $\delta$ induces a polar foliation of $S^n$ \cite[Theorem~4.2]{Thor}.
By \cite[Theorem~1.6]{Ly_resol}, the boundary $\D\Delta$ is precisely the set of singluar leaves. Hence, either  the claim holds or $\Delta$ has no boundary. 
In the latter case $\delta$ would be a regular Riemannian foliation with closed fibers.
Such a foliation is given by an $S^1$-action, an $S^3$-action or the Hopf fibration $S^{15}\to S^8$ \cite[Corollary~1.2]{LyWi_foliations}.
In all cases, the quotient would have the wrong curvature.
\qed

\bprop\label{prop_splitsub}
Let $\hat\delta:\hat\si\to\Delta$ be a transnormal submetry of a round $n$-sphere $\hat\si$. Suppose that there exists $(n+1)$
reflections $f_i$ along hyperspheres $\si_i\subset S^n$ with $\bigcap_{i=1}^{n+1}\si_i=\emptyset$ such that $\hat\delta$
is invariant under all $f_i$, $1\leq i\leq n+1$. Then either every geodesic in $\hat\si$  is piecewise $\hat\delta$-horizontal,
or there exists a disjoint decomposition $I^-\dot\cup I^+=\{1,\ldots,n+1\}$ such that $\hat\si=\si^-\circ\si^+$
where $\si^\pm=\bigcap_{i\in I^\pm}\si_i$ and every geodesic $\xi^-\xi^+$ with $\xi^\pm\in\si^\pm$ is piecewise $\hat\delta$-horizontal. 
\eprop

\proof
Let $\xi_1,\ldots,\xi_{n+1}\in \hat\si$ points such that $\xi_i\in\bigcap_{j\neq i}\si_j$.
If the links $\Si_{\xi_i}\hat\si$ are horizontal for $1\leq i\leq n+1$, then all links $\Si_\xi\hat\si$ are horizontal and every geodesic
is piecewise $\hat\delta$-horizontal since $\hat\delta$ is transnormal and the directions $\log_\xi(\xi_i)$
span the tangent space at $\xi$.
So let us assume that $\Si_{\xi_1}\hat\si$ is not horizontal and therefore splits into horizontal and vertical directions,
  $\Si_{\xi_1}\hat\si= H_{\xi_1}\circ V_{\xi_1}$. By assumption, this splitting is invariant under all reflections
	$f_i$, $2\leq i\leq n+1$. Set $v_i:=\log_{\xi_1}(\xi_i)$. It follows that there is a disjoint decomposition $I\dot\cup J=\{2,\ldots,n+1\}$
	such that the directions $v_i$ with $i\in I$ are vertical and the directions $v_j$ with $j\in J$ are horizontal.
	By assumption $I$ is not empty.
	Before we proceed, we need to make sure that vertical directions stay vertical.

\bslem
Let $\Pi$ be a fiber of $\delta$ and $x\in\Pi$ a point such that $T_x\Pi\subset T_x\si_i$ for some $i\in\{1,\ldots, n+1\}$.
Then the connected component $\Pi_0$ of $\Pi$ which contains the point $x$ has to lie in $\si_i$.
\eslem	

\proof
We may assume $i=1$.
By Proposition~\ref{prop_C11}, $\Pi$ is a $\mathcal{C}^{1,1}$ submanifold of $\hat\si$.
Let $\si\subset\si_1$ be a round sphere with $T_x\si=T_x\Pi$.
By the implicit function theorem, we can write $\Pi$ locally near $x$ as the graph of a $\mathcal{C}^{1}$
function $\varphi$ over $\si$. Since $\Pi$ is invariant under $f_1$, so is $\varphi$.
Therefore $\varphi$ takes values in $\si_1$. This shows that $\Pi$ lies in $\si_1$ locally near $x$.
Note that the invariance under $f_1$ implies that at an intersection point $y\in\Pi\cap\si_i$
the tangent space $T_y\Pi$ is either contained in $T_y\si_1$ or perpendicular to it.
Thus the condition that  a point $z$ of $\Pi$ lies in $\si_1$ and has tangent space $T_z\Pi$
contained in $T_z\si_1$ is open and closed. This implies the claim. 
\qed

A direct consequence is that every segment $\xi_k\xi_l$ is either vertical or horizontal.
Now denote by $\Pi$ the fiber of $\hat\delta$ through the point $\xi_1$.
Then the directions $v_i$, $i\in I$ form a basis of the tangent space $T_{\xi_1}\Pi$.
The sublemma implies that $\Pi$ is contained in the round sphere $\tilde\si:=\bigcap_{j\in J}\si_j$. 
Since $\dim(\Pi)=|I|=\dim(\tilde\si)$, it follows that $\Pi=\tilde\si$.
Thus, for $k,l\in I$, the segment $\xi_k\xi_l$ is vertical. 
We claim that for $k\in I$ and $l\in J$ the segment $\xi_k\xi_l$ is horizontal.
Recall that $\xi_1\xi_k$ is vertical and $\xi_1\xi_l$ is horizontal.
Thus if the segment $\xi_k\xi_l$ would be vertical, the triangle $\triangle(\xi_1,\xi_k,\xi_l)$
would entirely lie in $\Pi$. But then $\xi_1\xi_l$ cannot be horizontal.
Now set $I^+:=I\cup\{\xi_1\}$ and $I^-:=J$. Let $\xi_k\xi_l$ be a segment with $k\in I^+$
and $l\in I^-$. Since $I$ is non-empty, we can pick another point $\xi_j$ with $j\in I^+\setminus\{k\}$.
Then, the segments $\xi_k\xi_l$ and $\xi_j\xi_l$ are horizontal and the segment $\xi_j\xi_k$ is vertical.
Since the splitting into horizontal and vertical directions is orthogonal, the triangle $\triangle(\xi_j,\xi_k,\xi_l)$
has two right angles at $\xi_j$ and $\xi_k$. In particular $|\xi_j,\xi_l|=|\xi_k,\xi_l|=\frac{\pi}{2}$. 
This proves the claim.
\qed

\subsection{Spaces with an upper curvature bound}

For $\kappa \in \R $, let $D_{\kappa}  \in (0,\infty] $ be the diameter of the  complete, simply connected  surface $M^2_{\kappa}$
of constant curvature $\kappa$. A complete  metric space is called a \emph{CAT($\kappa$) space}
if any pair of its points with distance less than $D_{\kappa}$ is connected by a geodesic and if
 all triangles
with perimeter less than $2D_{\kappa}$
are not thicker than
 the \emph{comparison triangle} in $M_{\kappa} ^2$. In particular, geodesics between points of distance less than $D_{\kappa}$
are unique. 

For any CAT($\kappa$) space  $X$, 
 and points $x\neq y$ at distance less than $D_{\kappa}$, we denote the unique geodesic
 between $x$ and $y$ by $xy$.   For $y,z \neq x$, the angle at $x$ between $xy$ and $xz$
will be denoted by $\angle_x(y,z)$.
The \emph{space of directions} or {\em link}
 at a point  $x\in X$ is denoted by $\Si_x X$, it is a CAT(1) space when equipped with the angle metric.
Its elements are called {\em directions}.
A direction $v$ is called {\em genuine},
if there exists a geodesic starting in direction $v$.

\subsection{Dimension}
 A natural notion of dimension $\dim (X)$ for a CAT($\ka$) space $X$ was
  introduced by
Kleiner in \cite{Kleiner}, originally referred to as {\em geometric dimension}.
It vanishes precisely when the space is discrete.
In general, it is defined inductively:
\[\dim (X)= \sup _{x\in X} \{ \dim (\Sigma _xX) +1 \}.\]
For instance, in a 1-dimensional CAT($\ka$) space every link is discrete. 
The geometric dimension coincides with the supremum of 
topological dimensions\footnote{Here topological dimension corresponds to Lebesgue covering dimension.} of compact subsets in $X$ \cite{Kleiner}. 
If the dimension of $X$ is finite, then $\dim(X)$ agrees with the largest number $k$
such that $X$ admits a bilipschitz embedding of an open subset in $\R^k$ \cite{Kleiner}.
The dimension of a locally compact and geodesically complete space is finite and agrees with the topological dimension 
as well as the Hausdorff dimension \cite{OT_cba, LN_gcba}.  

\subsection{CAT(1) spaces}

CAT(1) spaces play a particular role in our investigations as they appear as Tits boundaries and links of CAT(0) spaces.
Recall that the spherical join $Z_1\circ Z_2$ of two CAT(1) spaces $Z_1$ and $Z_2$ is a CAT(1) space of diameter $\pi$.

 A subset $C$ in a CAT(1) space $Z$ is called {\em $\pi$-convex}
if for any pair of points $x, y\in Z$ at distance less than $\pi$ the unique geodesic $xy$ is contained in $C$.
If $C$ is closed, then it is CAT(1) with respect to the induced metric.
For instance, a ball of radius at most $\frac{\pi}{2}$ is $\pi$-convex.
Let $C\subset Z$ be a closed convex subset with radius $\rad(C)\geq \pi$. Then we define the set of {\em poles} for $C$ by
\[\pol(C):=\{\eta\in Z|\ d(\eta,\cdot)|_{C}\equiv\frac{\pi}{2}\}.\]
If $\diam(C)>\pi$, then $C$ has no pole. 
The convex hull of $C$ and $\pol(C)$ is canonically isometric to $C\circ\pol(C)$.

A subset  in a CAT(1) space is called {\em spherical}, if it 
embeds isometrically into a round sphere. By the Lune Lemma \cite[Lemma~2.5]{BB_diam}, a geodesic bigon in a CAT(1)
space is a spherical subset. Two points in a CAT(1) space are called {\em antipodes}, if their distance is at least $\pi$.
A point in a CAT(1) space is called {\em $k$-spherical}, if it has a neighborhood isometric to an open set in the round $k$-sphere.

\subsection{CAT(0) spaces}

The {\em ideal boundary} of a CAT(0) space $X$, equipped with the cone topology, is denoted by $\geo X$. 
If $X$ is locally compact, then $\geo X$ is compact.   
The {\em Tits boundary} of $X$ is denoted by $\tits X$, it is the ideal boundary equipped with 
the Tits metric $|\cdot,\cdot|_T$ -- the intrinsic metric associated to the {\em Tits angle}. 
 The Tits boundary of a CAT(0) space is a CAT(1) space.
 For any point $x\in X$ there is a  natural 1-Lipschitz  {\em logarithm map}
\[ \log_x:\tits X\to\Si_x X.\]
It satisfies the following rigidity. If $\log_x$ is isometric on a geodesic $\si$ in $\tits X$, then 
the geodsic cone from $x$ to $\si$ is isometric to a flat sector \cite[Flat Sector Lemma~2.3.4]{KleinerLeeb}.
 
If $X_1$ and $X_2$ are CAT(0) spaces, then their metric product $X_1\times X_2$ is again a CAT(0) space.
We have $\tits (X_1\times X_2)=\tits X_1\circ \tits X_2$ and $\Si_{(x_1,x_2)}(X_1\times X_2)=\Si_{x_1} X_1\circ \Si_{x_2} X_2$.  
If $X$ is a geodesically complete CAT(0) space, then any join decomposition of $\tits X$ is induced by a metric product decomposition of $X$ \cite[Proposition~2.3.7]{KleinerLeeb}.

A {\em $n$-flat} $F$ in a CAT(0) space $X$ is a closed convex subset isometric to $\R^n$.
Its Tits boundary $\tits F\subset\tits X$ is a round $(n-1)$-sphere.
On the other hand, if $X$ is locally compact and $\si\subset \tits X$ is a round  $(n-1)$-sphere, then either there exists an
$n$-flat $F\subset X$ with $\tits F=\si$, or there exists a round $n$-hemisphere $\tau^+\subset\tits X$
with $\si=\D\tau^+$ \cite[Proposition~2.1]{Leeb}. Consequently, if $\tits X$ is $(n-1)$-dimensional, then any round $(n-1)$-sphere in $\tits X$
is the Tits boundary of some $n$-flat in $X$.  A {\em flat ($n$-dimensional) half-space} $H\subset X$ is a closed convex subset isometric to a Euclidean half-space $\R_+^n$.

\subsection{Basics on ideal boundaries}

We record two facts which relate the cone topology and the Tits metric.

\blem\label{lem_angle_to_set}
Let $X$ be a locally compact  CAT(0) space.
Let $A\subset \geo X$ be a closed subset and $\eta$ a point in $\geo X\setminus A$
represented by a ray $\rho$. Then $\lim\limits_{t\to\infty}\angle_{\rho(t)}(\eta, A)=|\eta, A|$ 
holds.
\elem

\proof
Since $\log_x:\tits X\to\Si_x X$ is 1-Lipschitz for every point $x\in X$, we have $\limsup\limits_{t\to\infty}\angle_{\rho(t)}(\eta, A)\leq|\eta, A|$.
So if the statement fails, then we find $\epsilon>0$, $t_i\to\infty$ and $\xi_i\in A$ with $\xi_i\to\xi\in A$ such that for all $i\in\N$ holds
\[\angle_{\rho(t_i)}(\eta, \xi_i)\leq|\eta,A|-\epsilon.\]

But for any $\delta>0$ we find $i_0\in\N$ such that $|\eta,\xi|\leq\angle_{\rho(t_i)}(\eta,\xi)+\delta$ for
$i\geq i_0$. Choose $j_0\geq i_0$ such that $\angle_{\rho(t_{i_0})}(\eta,\xi)\leq \angle_{\rho(t_{i_0})}(\eta,\xi_j)+\delta$ for $j\geq j_0$.
Together, we obtain 
\[
|\eta,\xi|\leq\angle_{\rho(t_{i_0})}(\eta,\xi_j)+2\delta\leq\angle_{\rho(t_j)}(\eta,\xi_j)+2\delta\leq|\eta,A|-\epsilon+2\delta.
\]
This leads to a contradiction if we choose $\delta<\frac{\eps}{2}$.
\qed

\blem\label{lem_mindist}
Let $X$ be a locally compact CAT(0) space.
Let $A^\pm$ be closed subsets of $\geo X$. Then there exist
 points $\xi^\pm\in A^\pm$ with 
\[|\xi^-,\xi^+|=|A^-, A^+|.\]
\elem

\proof
If there exists $\eta\in A^-\cap A^+$, we can take $\xi^-=\xi^+=\eta$. Otherwise,
there exist sequences of points $\xi^\pm_i\in A^\pm,$ such that $|\xi_i^-,\xi_i^+|\to|A^-, A^+|$.
Since $\geo X$ is compact, we may assume that $\xi^\pm_i\to\xi^\pm$ with respect to the cone topology.
By lower semicontinuity of the Tits metric, we conclude $|\xi^-,\xi^+|\leq|A^-, A^+|$.
By assumption, $A^\pm$ is closed and therefore $\xi^\pm\in A^\pm$.
\qed

\subsection{Dynamics at infinity}\label{subsec_dyn}
Recall the following construction from \cite{GuSw_trans}.
Let $G$ be a discrete group acting on a compact Hausdorff space $Z$.
Identify the Stone-Čech compactification $\beta G$ with the set of ultrafilters on 
$G$. For every $\om\in\beta G$ define 
\[T^\om:Z\to Z;\ z\mapsto\wlim_g gz. \] 
Then for fixed $z_0\in Z$ the map $\beta G\to Z$ which sends $\om$
to $T^\om z_0$ is continuous. The family of operators $\{T^\om\}_{\om\in\beta G}$
is closed under composition and the inverse map $g\mapsto g^{-1}$
extends to a continuous involution $S:\beta G\to\beta G$.

Now let $X$ be a locally compact CAT(0) space and $\Ga\acts X$ a geometric action.
Then $\bar X=X\cup\geo X$ is a compact Hausdorff space and the above construction
applies. By the semi-continuity of the Tits metric, every operator $T^\om$ is a 1-Lipschitz map $\tits X\to\tits X$ \cite{GuSw_trans}.

For $\om\in\beta\Ga\setminus\Ga$ and $x\in X$ we denote by $\om^+:=T^\om(x)$
and $\om^-:=T^{S\om}(x)$ the {\em attracting} and {\em repelling} point, respectively. The definition does not depend on the choice of $x\in X$.

The key tool for studying the dynamics of $\Ga$ on $\geo X$
is the following theorem due to Papasoglu and Swenson.

\bthm[{$\pi$-convergence \cite[Lemma~19]{PaSw_bound}}]
Let $\Ga$ be a group acting properly discontinuously on a locally compact CAT(0)
space $X$, and let $\om\in\beta\Ga\setminus\Ga$, then
\[|\om^-,\xi|+|T^\om(\xi),\om^+|\leq\pi\]
holds for every $\xi\in\tits X$.
\ethm

\proof
It is enough to show the following. Suppose that $(\ga_k)$ is a sequence in $\Ga$
with $\ga_k^{\pm 1} x\to\om^\pm$ for some (and then any) point $x\in X$ and such that $\ga_k\xi\to\eta$ in $\geo X$.
Then for any point $y\in X$ we have
\[\angle_y(\ga_k^{-1}x,\xi)=\angle_{\ga_k y}(x,\ga_k\xi)\leq\pi-\angle_{x}(\ga_k y,\ga_k\xi).\]
Now choose $y_k\in y\xi\setminus\{y\}$ with $y_k\to y$ and $x_k\in x\ga_k\xi\setminus\{x\}$ with $x_k\to x$.
Passing to the limit $k\to\infty$ and using \cite[Lemma~2.1.5]{KleinerLeeb}, we conclude
\[\angle_y(\om^-,\xi)\leq\pi-\angle_{x}(\om^+,\eta).\]
Taking the supremum over $y\in X$ and the infimum over $x\in X$ yields the result.
\qed
\medskip

An element $\om\in\beta\Ga$ is said to {\em pull} from a point $\xi\in\geo X$, if 
the following holds. Let $x\in X$ and $\rho:[0,\infty)\to X$ be a geodesic ray asymptotic to $\xi$.
Then there exists a constant $r>0$ such that the set of all $\ga\in\Ga$ with $\ga x\in N_r(\rho')$ lies in $S\om$
for every subray $\rho'\subset\rho$. In particular, the repelling point is $\om^-=\xi$.
Here is the basic example. Suppose $(\ga_k)$ is a sequence in $\Ga$ such that $(\ga_k\rho(k))$
is bounded.
If we push forward a non-principal ultrafilter on $\N$  to $\Ga$
via $(\ga_k)$, then we obtain an element $\om$ which pulls from $\xi$.
In particular, if $\Ga$ acts geometrically then we can pull from every point in $\geo X$.
We will make use of the following.

\blem[{\cite[Lemma~3.19]{GuSw_trans}}]\label{lem_pullflat}
Let $\si\subset\tits X$ be a round $n$-sphere which bounds a flat $F\subset X$.
If $\om\in\beta\Ga$ pulls from a point $\xi\in\si$, then $T^\om\si$
is a round $n$-sphere which also bounds a flat.
\elem  

A subset $A\subset\tits X$ is called {\em incompressible}, if $T^\om$ maps $A$ isometrically
for every $\om\in\beta\Ga$. Note that a maximal incompressible set is automatically Tits-closed and convex.

\subsection{Morse flats}

Morse quasiflats in metric spaces were introduced and studied in \cite{HKS_I,HKS_II}. They simultaneously generalize Morse quasi-geodesics as well as top-dimensional quasiflats.
Their significance lies in the fact that they are quasi-isometry invariant and  display nice asymptotic behavior, namely they 
have unique tangent cones at infinity. This makes them a strong tool in asymptotic geometry, for instance in the study of quasi-isometric rigidity of infinite groups.
Here we are only concerned with the very particular case of periodic Morse flats. Recall that a flat $F$ in a metric space $X$ is {\em periodic}, if its stabilizer in the full isometry group of $X$ contains a subgroup which acts geometrically on $F$.
A periodic $n$-flat $F$ is called {\em Morse}, if it does not bound a flat $(n+1)$-half-space.
Periodic Morse $1$-flats are precisely Ballmann's rank 1 axes \cite[Chapter~III.3]{ballmannbook}. If  
$\ga$ is a periodic geodesic in a CAT(0) space $X$, then $\ga$ is Morse if and only if (one and then) both ideal endpoints of
$\ga$ are isolated in the Tits boundary. In higher dimensions, we have the following, which follows immediately from 
\cite[Corollary~11.5]{HKS_I} and \cite[Proposition~6.14]{HKS_II}.

\bprop\label{prop_morseinj}
Let $X$ be a locally compact CAT(0) space and $F\subset X$ a periodic $(n+1)$-flat with $\si:=\tits F$.
Then $F$ is Morse if and only if for every point $\xi\in\si$ the induced map
\[\tilde H_n(\si,\si\setminus\{\xi\},\Z)\to \tilde H_n(\tits X,\tits X\setminus\{\xi\},\Z)\]
is injective. 
\eprop 

An important consequence is that any ideal boundary point has an antipode in the Tits
boundary of a periodic Morse flat. 

\bcor\label{cor_morseantipodes}
Let $X$ be a locally compact CAT(0) space and $F\subset X$ a periodic Morse flat with $\si:=\tits F$.
Then any point in $\tits X$ has an antipode in $\si$. Moreover, any geodesic in $\tits X$ which ends in $\si$
can be extended inside $\si$.
\ecor

\proof
Let $\ga$ be a geodesic with endpoint $\xi\in\si$. Denote by $v\in\Si_\xi\tits X$ the incoming direction.
If $\Si_\xi\si\subset B_\pi(v)$, then $\Si_\xi\si$ is contractible inside $\Si_\xi\tits X$ since every ball of radius less than $\pi$
is contractible in a CAT(1) space. By \cite[Theorem~A]{Kramer}, the logarithm map provides a homotopy equivalence between a small enough punctured ball around
$\xi$ and $\Si_\xi\tits X$. In particular, the map $\tilde H_n(\si,\si\setminus\{\xi\},\Z)\to \tilde H_n(\tits X,\tits X\setminus\{\xi\},\Z)$
is trivial. Contradiction. Thus we do find an antipode of $v$ in $\Si_\xi\si$.
\qed
\medskip

The following is similar to \cite[Corollary~2.5]{BMS_affine}.

\bcor\label{cor_morsecyclesdense}
Let $X$ be a locally compact CAT(0) space with a geometric group action $\Ga\acts X$ and let $F\subset X$ be a periodic Morse flat with $\si:=\tits F$.
Then the orbit $\Ga\si$ is dense in $\geo X$.
\ecor

\proof
Let $\xi$ be a point in $\geo X$. Choose an element $\om\in\beta\Ga$ with repelling point $\xi=\om^-$.
By Corollary~\ref{cor_morseantipodes}, the attracting point $\om^+$ has an antipode $\eta$ in $\si$.
By $\pi$-convergence, $T^\om(\eta)=\xi$. Since $\xi$ has a countable neighborhood basis in $\bar X$
we obtain a sequence $(\ga_k)$ in $\Ga$ with $\ga_k\eta\to\xi$.  
\qed

\begin{rem}
In view of \cite[Proposition~6.14]{HKS_II} it is clear that the conclusions of both corollaries remain true for Tits boundaries of Morse (quasi)-flats,
no periodicity required.
\end{rem}

\blem\label{lem_morseretract}
Let $X$ be a locally compact CAT(0) space with a geometric group action $\Ga\acts X$.
Suppose that $F$ is a $\Ga$-periodic Morse flat with $\si=\tits F$.
Then there exists $\om\in\beta\Ga$ such that $T^\om$ is a 1-Lipschitz retraction $\tits X\to\si$.
\elem

\proof
By Theorem~\ref{thm_regular}, the set of regular points in $\si$ is open and dense.
Choose an axial isometry $\ga$ in the stabilizer of $F$ with a regular axis $c\subset F$.
Push forward a non-principal ultrafilter from $\N$ to $\Ga$ via $k\mapsto \ga^k$.
Then $\om^\pm=c(\pm\infty)$ and $\om$ pulls from $c(-\infty)$.

\qed

\section{Symmetric subsets}\label{sec_sym}

\subsection{Antipodes and visual antipodes}

Let us refine the concept of antipodal points to take into account the 
space $X$.

\bdfn
Let $X$ be a  CAT(0) space.
Two points $\xi,\hat\xi\in\geo X$ are called {\em antipodal}, if $|\xi,\hat\xi|\geq\pi$. They are called
{\em visually antipodal}, if there is a complete geodesic $c$ in $X$ with $\geo c=\{\xi,\hat\xi\}$. For a subset
$A\subset \geo X$ we denote the set of (visually) antipodal points of $A$ by $\Ant_{(vis)}(A)$. Inductively,
we define $\Ant^{j+1}_{(vis)}(A)=\Ant_{(vis)}(\Ant^{j}_{(vis)}(A))$. 
\edfn

\blem\label{lem_dense}
Let $X$ be a locally compact and geodesically complete CAT(0) space.
 For any subset $A\subset\geo X$ holds:
\begin{enumerate}
 \item[(i)] $\Ant_{vis}(A)$ is dense in $\Ant(A)$;
 \item[(ii)] If $A$ is dense in $A'$, then $\Ant(A)$ is dense in $\Ant(A')$.
\end{enumerate}
\elem

\proof (i)
Let $\eta$ be a point in $A$ and $\hat\eta\in\Ant(\eta)$. Represent $\hat\eta$ 
by a ray $\hat \rho$. For a sequence $t_i\to\infty$ extend the rays $\hat \rho(t_i)\eta$
to complete geodesics $c_i$, with ideal boundary points $c_i(-\infty)=\eta$
and $c_i(+\infty)=\hat\eta_i$, such that 
$\angle_{\hat \rho(t_i)}(\eta,\hat\eta_i)=\pi-\angle_{\hat \rho(t_i)}(\eta,\hat\eta)$.
Then $\hat\eta_i\to\hat\eta$, since for every point $p$ on $\hat\rho$ holds 
$\angle_{p}(\eta,\hat\eta_i)\leq\angle_{\hat \rho(t_i)}(\eta,\hat\eta_i)=
\pi-\angle_{\hat \rho(t_i)}(\eta,\hat\eta)\to 0$.
\medskip

(ii) 
Direct consequence of (i) since $\Ant_{vis}(A)$ is dense in $\Ant_{vis}(A')$.
\qed
\medskip

By induction, we obtain

\bcor\label{cor:dense}
Let $X$ be a locally compact and geodesically complete CAT(0) space.
Let $A$ be a subset of $\geo X$. Then $\Ant^j_{vis}(A)$ is dense in $\Ant^j(A)$
for every $j\in\N$.
\ecor

\bdfn
Let $X$ be a CAT(0) space.
A subset $A$ in $\geo X$ is called {\em (visually) symmetric}, if 
\[\Ant_{(vis)}(A)\subset A.\]
A visually symmetric subset is also called {\em involutive} \cite{E_sym}.
\edfn

\bcor\label{cor_vis}
Let $X$ be a locally compact and geodesically complete CAT(0) space.
A closed subset $A\subset\geo X$ is symmetric if and only if it is visually symmetric.
\ecor

\proof
Any symmetric subset is visually symmetric. By Lemma~\ref{lem_dense}, for a closed visually symmetric subset $A\subset\geo X$ holds
\[
\Ant(A)\subset\overline{\Ant(A)}=\overline{\Ant_{vis}(A)}\subset A.
\]
\qed

\bcor\label{cor:closuresym}
Let $X$ be a locally compact and geodesically complete CAT(0) space.
If a subset $A\subset\geo X$ is visually symmetric, then its closure $\overline{A}$ is symmetric.
\ecor

\proof
$\Ant_{vis}(A)$ is dense in $\Ant(\overline A)$ by Lemma~\ref{lem_dense}. Therefore
\[
\Ant(\overline A)\subset\overline{\Ant(\overline A)}=\overline{\Ant_{vis}(A)}\subset\overline A.
\] 
\qed

\blem\label{lem_bound}
Let $X$ be a locally compact and geodesically complete CAT(0) space.
Let $A$ be a closed symmetric subset of $\geo X$ and $\eta$ a point in the complement $\geo X\setminus A$. Then  for every $\xi\in A$ holds 
\[|\eta,\xi|\leq\pi-|\eta,A|.\]
In particular, $|\eta,A|\leq\frac{\pi}{2}$. 
\elem

\begin{center}
\includegraphics[scale=0.35,trim={0cm 4cm 0cm 5cm},clip]{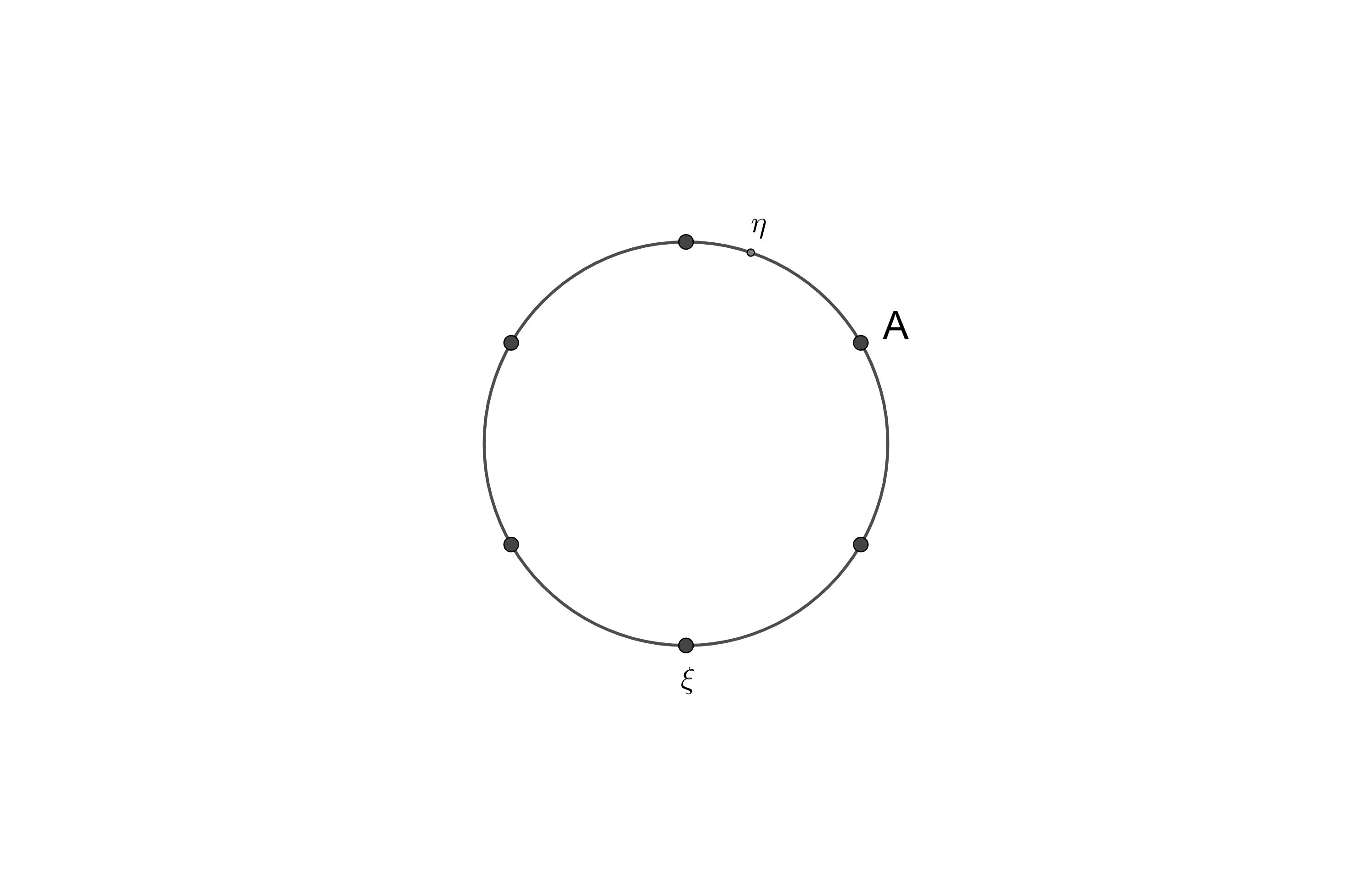}
\end{center}

\proof
Let $\xi$ be a point in $A$ and $\rho$ a ray representing $\eta$. For every $t\geq 0$
extend the ray $\rho(t)\xi$ to a complete geodesic with ideal endpoints $\hat\xi_t$ and $\xi$
such that $\angle_{\rho(t)}(\hat\xi_t,\eta)+\angle_{\rho(t)}(\eta,\xi)=\pi$. Since $A$ is symmetric, $\hat\xi_t\in A$.
Then for $t\geq 0$ holds
\[\angle_{\rho(t)}(\xi,\eta)=\pi-\angle_{\rho(t)}(\eta,\hat\xi_t)\leq \pi-\angle_{\rho(t)}(\eta,A).\]
For $t\to\infty$, the left hand side converges to $|\xi,\eta|$. By Lemma~\ref{lem_angle_to_set},
the right hand side converges to $\pi-|A,\eta|$. The claim follows.
\qed

\bcor\label{cor_imagedist}
Let $X$ be a locally compact and geodesically complete CAT(0) space.
Let $A$ be a closed symmetric subset of $\geo X$ and let $f:\tits X\to Z$ be a surjective 1-Lipschitz map to a geodesically complete 
CAT(1) space $Z$. Then for every point $\eta\in\geo X$ holds $|f(\eta),f(A)|=|\eta,A|$. 
Moreover, if $|\eta, A|=|\eta,\xi|$ holds for some point  $\xi\in A$, then $|f(\eta),f(A)|=|f(\eta),f(\xi)|=|\eta,\xi|$.
\ecor

\proof
Suppose for contradiction that there is a point $\xi\in A$ with $|f(\eta),f(\xi)|<|\eta,A|$.
Since $f$ is surjective and 1-Lipschitz, and $Z$ is geodesically complete, we find an antipode $\hat\xi\in A$ of $\xi$
such that $|f(\hat \xi),f(\eta)|+|f(\eta),f(\xi)|=\pi$. From Lemma~\ref{lem_bound}, we conclude
\[|\eta,A|\leq\pi-|\hat\xi,\eta|\leq \pi-|f(\hat\xi),f(\eta)|=|f(\eta),f(\xi)|<|\eta,A|.\]
Contradiction. The supplement follows from
\[|\eta,\xi|=|\eta,A|=|f(\eta),f(A)|\leq|f(\eta),f(\xi)|\leq|\eta,\xi|.\]
\qed

If we apply this to $f=\log_o$, we obtain:

\bcor\label{cor_prefact}
Let $X$ be a locally compact and geodesically complete CAT(0) space.
Let $A\subset \geo X$ be a closed symmetric subset and $\eta$ a point in $\geo X$.
Then for every point $o\in  X$ holds $\angle_o(\eta,A)=|\eta,A|$.
If $|\eta, A|=|\eta,\xi|$ holds for some point  $\xi\in A$, then $\angle_o(\eta,A)=\angle_o(\eta,\xi)=|\eta,\xi|$.
\ecor

Another immediate consequence is

\bcor\label{cor_imagedistsets}
Let $X$ be a locally compact and geodesically complete CAT(0) space.
Let $f:\tits X\to Z$ be a surjective 1-Lipschitz map to a geodesically complete 
CAT(1) space $Z$.
If $A,B\subset\geo X$ are equidistant closed symmetric subsets, then so are $f(A),f(B)\subset Z$, and 
$|A,B|=|f(A),f(B)|$ holds. Moreover, if $|A,B|=|\xi,B|$ holds for some point $\xi\in A$, then 
$|f(A),f(B)|=|f(\xi),f(B)|=|\xi,B|$.
\ecor

\subsection{Almost near points}\label{sec_alnear}

\bdfn
Let $Z$ be a CAT(1) space and $\xi,\eta$ two points in $Z$.
Then $\eta$ is {\em almost near} to $\xi$, if for every antipode $\hat\xi$ of $\xi$ holds
\[|\xi,\eta|+|\eta,\hat\xi|=\pi.\] 
\edfn

\begin{rem}
Note that if an antipode $\hat\xi$ of $\xi$ is almost near to $\xi$, then $\hat\xi$ is the only antipode 
of $\xi$. 
Further, if a CAT(1) space $Z$ where every point has an antipode contains a pair of points $\xi,\eta$ such that $\eta$ is almost near to $\xi$,
then the geodesic $\xi\eta$ extends up to distance $\pi$ and the $\pi$-ball around $\xi$ covers all of $Z$.
This applies in particular to Tits boundaries of geodesically complete CAT(0) spaces. 
\end{rem}

\bcor\label{cor_aln}
Let $X$ be a locally compact and geodesically complete CAT(0) space.
Let $A$ be a closed symmetric subset of $\geo X$ and $\eta$ a point in the complement $\geo X\setminus A$.
If $\xi\in A$ realizes the distance between $\eta$ and $A$, $|\eta, A|=|\eta,\xi|$, then $\eta$ is
almost near to $\xi$. 
\ecor

\proof
Since $A$ is symmetric, every antipode $\hat\xi$ of $\xi$ lies in $A$. Hence, Lemma~\ref{lem_bound} implies
\[\pi\leq|\xi,\eta|+|\eta,\hat\xi|\leq |\xi,\eta|+\pi-|A,\eta|=\pi.\]
\qed

\blem\label{lem_aln_iso}
Let $X$ be a locally compact and geodesically complete CAT(0) space.
Let $\xi,\eta\in\geo X$ be points such that $\eta$ is almost near to $\xi$.
Then for any point $o\in X$ we have 
\[\angle_o(\xi,\eta)=|\xi,\eta|\]
and the point $\log_o(\eta)$ is almost near to $\log_o(\xi)$ in $\Si_o X$.
\elem

\proof
For any antipode $\hat v$ of $\log_o(\xi)$ we can extend the ray $o\xi$ in the direction $\hat v$ 
to a complete geodesic $c$ with ideal boundary $\geo c=\{\xi,\hat\xi\}$.
Since $\eta$ is almost near to $\xi$, there exists a geodesic from $\xi$ to $\hat\xi$
 through the point $\eta$. Hence $c$ bounds a flat half-plane $H$ with $\eta\in\geo H$ \cite[Flat Sector Lemma~2.3.4]{KleinerLeeb}.
In particular, 
\[\pi\leq\angle_o(\xi,\eta)+\angle_o(\eta,\hat v)\leq|\xi,\eta|+|\eta,\hat\xi|=\pi\]
and the claim follows.
\qed

\bprop\label{prop_sph}
Let $\xi,\eta$ be points in $\geo X$ such that $\eta$ is almost near to $\xi$.
Let  $\zeta\in\geo X$ be a point such that the geodesic triangle $\triangle(\xi,\eta,\zeta)$ has perimeter less than $2\pi$. 
Then  $\triangle(\xi,\eta,\zeta)$ is spherical.
\eprop

\begin{center}
\includegraphics[scale=0.3,trim={-3cm 7cm 0cm 1.5cm},clip]{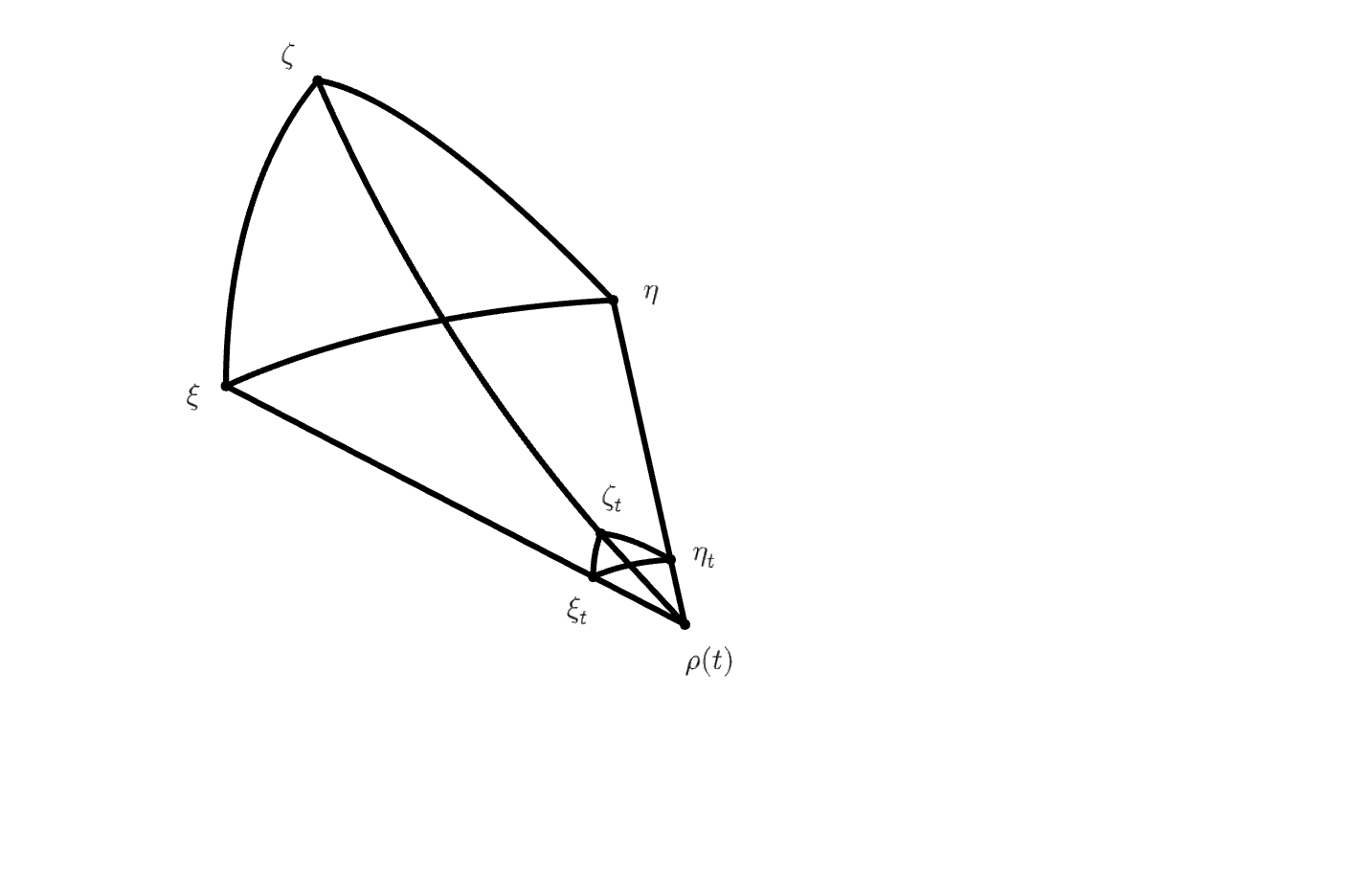}
\end{center}

\proof
By Reshetnyak majorization, we find a 1-Lipschitz map $\varphi:\tilde\triangle\to \triangle(\xi,\eta,\zeta)$
where $\tilde\triangle$ denotes the comparison triangle and $\varphi$ maps every side isometrically.
Let $\rho$ be a geodesic ray representing $\zeta$. For every $t\geq 0$ we consider the 1-Lipschitz map 
\[f_t=\log_{\rho(t)}\circ\varphi:\tilde\triangle\to\Si_{\rho(t)}X\]
and an ultralimit with respect to a non-principal ultrafilter $\om$: 
\[f_\om=\wlim f_k:\tilde\triangle\to\Si_\om\] 
where $\Si_\om=\wlim\Si_{\rho(k)}X$.
For a point in $\tits X$ we use a subscript $t$ to denote its image under $\log_{\rho(t)}$. For every $t\geq 0$, the link $\Si_{\rho(t)} X$ is geodesically complete
and, By Lemma~\ref{lem_aln_iso}, $\eta_t$ is almost near to $\xi_t$. Hence the triangle $\triangle(\xi_t,\eta_t,\zeta_t)$ is spherical \cite[Proposition~6.1]{Ly_rigidity}. Thus the ultralimit $\triangle_\om=\wlim \triangle(\xi_t,\eta_t,\zeta_t)$ is spherical as well. 
By Lemma~\ref{lem_aln_iso}, $f_t$ maps $\xi\eta$ isometrically, hence so does $f_\om$. 
Moreover, since 
\[|\zeta,\eta|=\lim\limits_{t\to\infty}\angle_{\rho(t)}(\zeta,\eta)\ \text{ and }\ |\zeta,\xi|=\lim\limits_{t\to\infty}\angle_{\rho(t)}(\zeta,\xi),\]
$f_\om$ maps all three sides of $\triangle(\xi,\eta,\zeta)$ isometrically. In particular, $f_\om$ provides an isometry $\tilde\triangle\to\triangle_\om$.
Since $f_\om$ factors through $\triangle(\xi,\eta,\zeta)$, we obtain the claim.
\qed

\blem\label{lem_nobranch}
Let $\xi,\eta$ be points in $\geo X$ such that $\eta$ is almost near to $\xi$.
Then no geodesic from $\xi$ to $\eta$ contains a branch point in the interior.
\elem

\proof
Suppose that $\beta$ is branch point in the interior of a geodesic $\ga$ from $\xi$ to $\eta$.
Then there exist points $\eta^-$ in the interior of the segment $\beta\eta$ and a point $\eta^+$
not on $\ga$, such that the geodesics $\xi\eta^-$ and $\xi\eta^+$ intersect in the segment $\xi\beta$;
and the perimeter of the triangle $\triangle(\xi,\eta^-,\eta^+)$ is less than $2\pi$.
But $\eta^-$ is almost near to $\xi$. Thus, by Proposition~\ref{prop_sph}, the triangle $\triangle(\xi,\eta^-,\eta^+)$
is spherical. Contradiction.
\qed

\blem\label{lem_2flat}
Let $X$ be a locally compact and geodesically complete CAT(0) space.
Let $A$ be a closed symmetric subset of $\geo X$ and $\eta$ a point in the complement $\geo X\setminus A$. 
Let $\xi\in A$ be a point which realizes the distance from $\eta$ to $A$, $|\eta, A|=|\eta, \xi|$. Suppose that $c$ is a 
complete geodesic in $X$, asymptotic to $\eta$. Then $c$ lies in
a 2-flat $F\subset X$ which contains $\xi$ in its boundary circle $\geo F$. In particular, for every $\hat\eta\in\Ant_{vis}(\eta)$ the segment $\xi\eta$
lies in a simple closed geodesic containing $\hat\eta$.
\elem

\proof
For every $t\in\R$ we extend $c(t)\xi$ to a complete geodesic $l_t$ with ideal 
endpoints $l_t(+\infty)=\xi$ and $l_t(-\infty)=\hat\xi_t$. By Corollary~\ref{cor_aln}, $\eta$ is almost near to $\xi$ and therefore
there exists a geodesic from $\xi$ to $\hat\xi_t$
 through the point $\eta$. Hence,
 $l_t$ bounds a flat half-plane $H_t$ whose boundary 
half-circle $\geo H_t$ contains $\eta$. For a sequence $t_i\to-\infty$ the
points $\hat\xi_{t_i}$ subconverge to a point $\hat\xi\in A$ and the half-planes
$H_{t_i}$ converge to a 2-flat $F$. Clearly $F$ contains $c$
and $\xi\in\geo F$.
\qed

%

\section{A submetry at infinity}\label{sec_sub}

\subsection{Minimal closed symmetric sets}

The goal of this section is to show that the existence of a non-trivial, proper, closed, symmetric subset $A\subset\geo X$
leads to a non-trivial submetry $\tits X\to \Delta$ (Proposition~\ref{prop_submetry}).

\bdfn
Let $X$ be a CAT(0) space.
 For a subset $U\subset\geo X$ we denote $A^{(vis)}_U$ the minimal (visually) symmetric subset
containing $U$. Moreover, for a point $\xi\in\geo X$ we denote by $B_\xi$ the closure of $A^{vis}_\xi$. By
Corollary~\ref{cor:closuresym}, the set $B_\xi$ is symmetric.

\noindent(Note: $A^{(vis)}_U=U\cup\Ant^{(vis)}(U)\cup(\Ant^{(vis)})^2(U)\cup\ldots$) 
\edfn

\blem\label{lem_const}
Let $X$ be a locally compact and geodesically complete CAT(0) space.
Let $A$ be a closed symmetric subset of $\geo X$ and $\eta$ a point in $\geo X$. 
Then $|\cdot,A|$ is constant on $B_\eta$. In particular, $|B_\eta,A|=|\eta,A|$ holds. 
\elem

\proof
Lemma~\ref{lem_2flat} implies
that for $\hat\eta\in\Ant^{vis}(\eta)$ there exists $\hat\xi\in A$ with 
$|\hat\eta,\hat\xi|=|\eta,A|$. Hence $|\hat\eta,A|\leq|\eta,A|$ holds, and by symmetry we have
$|\hat\eta,A|=|\eta,A|$. Therefore, $|\cdot,A|$ is constant on $A^{vis}_\eta$. 
By Lemma~\ref{lem_mindist}, we have $|B_\eta,A|=|\eta_\infty,\xi|$ for some points $\eta_\infty\in B_\eta$ and $\xi\in A$.
By Corollary~\ref{cor:dense}, we find a sequence of points $\eta_k\in A^{vis}_\eta$ such that $\eta_k\to \eta_\infty$.
We choose a base point $o\in X$.
Using Lemma~\ref{lem_aln_iso}, we see that for every $\eps>0$ and $k\in\N$ large enough, the following estimates hold
\begin{align*}
|B_\eta,A|&=|\eta_\infty,\xi|=\angle_o(\eta_\infty,\xi)\geq\angle_o(\eta_k,\xi)-\eps\\
&\geq\angle_o(\eta_k,A)-\eps=|\eta_k,A|-\eps=|\eta,A|-\eps.
\end{align*}
In the next to last step we have used Corollary~\ref{cor_prefact}.
Since $\eps>0$ was arbitrary, the claim follows.
\qed

\bcor\label{cor:presubmetry}
Let $X$ be a locally compact and geodesically complete CAT(0) space.
For $\xi,\eta\in\geo X$ holds $|B_\xi,\eta|=|B_\xi,B_\eta|=|\xi,B_\eta|$. In particular,
$\xi\in B_\eta$ and $\eta\in B_\xi$ are equivalent, and the set $B_\xi$ is the minimal closed symmetric set containing $\xi$.
\ecor

\begin{center}
\includegraphics[scale=0.4,trim={1cm 4cm 0cm 3cm},clip]{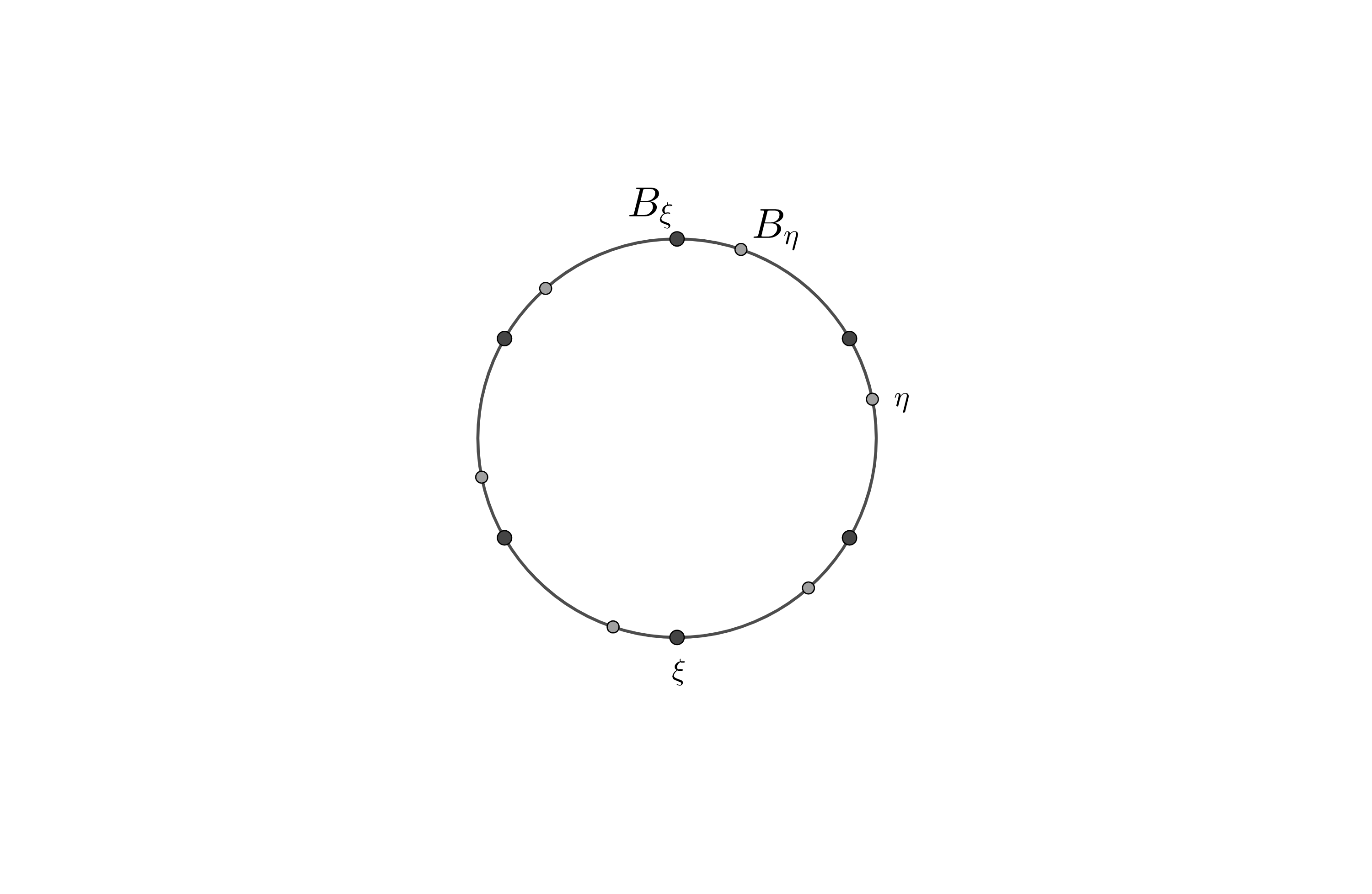}
\end{center}

\bprop\label{prop_submetry}
Let $X$ be a locally compact and geodesically complete CAT(0) space.
The sets $B_\xi$ constitute an   equidistant decomposition of $\tits X$ and therefore induce
a submetry 
\[\delta:\tits X\rightarrow \Delta.\]
If $\Delta$ is not a point, i.e. if $\geo X$ contains a non-empty, proper, closed, symmetric subset,  then the diameter of $\tits X$ is $\pi$.
\eprop

\proof
It follows from Corollary~\ref{cor:presubmetry} and
Lemma~\ref{lem_mindist} that $\delta$ is a submetry.
If the diameter of $\tits X$ is larger than $\pi$, then there exists a complete geodesic $c\subset X$
such that the ideal endpoints $\xi^\pm=c(\pm\infty)$ have distance larger than $\pi$.
Let $\eta\in\tits X$ be a point which is $\delta$-near to $\xi^+$. 
Then, by Lemma~\ref{lem_2flat}, there exists a 2-flat $F$ which contains $c$ and whose Tits boundary $\tits F$
contains the segment $\xi^+\eta$. Contradiction. Thus the diameter of $\tits X$ has to be $\pi$.
\qed
\medskip

%

An application of Corollary~\ref{cor_imagedistsets} with $f=\log_o$ leads to

\bcor\label{cor_fact}
Let $X$ be a locally compact and geodesically complete CAT(0) space.
For every point $o\in X$, the sets $\log_o(B_\xi)$ form an equidistant decomposition of $\Si_o X$
into closed symmetric subsets; and $\angle_o(\xi, B_\eta)=\angle_o(B_\xi, B_\eta)=|B_\xi, B_\eta|$ holds.
\ecor


\medskip

As a consequence, we obtain for every $o\in X$ a submetry 
\[\delta_{o}:\Si_o X\to\Delta.\]
Since the fibers of $\delta_{o}$ are closed symmetric subsets, the map has to factor
through the submetry $\mu_o:\Si_o X\to\Delta_o$ induced by the decomposition into minimal closed symmetric subsets.
So we obtain the following diagram where all maps but $\log_o$ are submetries.
\medskip

\adjustbox{scale=1.5,center}{%
\begin{tikzcd}[scale=3]
\tits X \arrow{d}{\delta} \arrow{r}{\log_o} & \Si_o X \arrow{d}{\mu_o}\arrow{dl}{\delta_{o}} \\
\Delta & \Delta_{o}\arrow{l}{\bar\delta_{o}}
\end{tikzcd}
}
\medskip

Since in our setting all links $\Si_o X$ are locally compact and geodesically complete,
\cite{Ly_rigidity} provides us with control on the submetries $\delta_{o}:\Si_o X\to\Delta$.
We will use this later on to derive properties of the submetry $\delta:\tits X\to\Delta$.
As a matter of fact, \cite{Ly_rigidity} implies that every space of directions $\Si_o X$ is
either a spherical join or a spherical building.

\subsection{Horizontal geodesics}\label{sec_hor}

In this section we begin our study of the submetry $\delta:\tits X\to\Delta$ provided by Proposition~\ref{prop_submetry}.
A key property is {\em horizontal geodesic completeness.}
Namely, from Corollary~\ref{cor_aln} we obtain:

\bprop\label{prop_delnearalnear}
If  a point $\eta\in\tits X$ is $\delta$-near to $\xi$, then it is almost near.
In particular, the geodesic $\xi\eta$ extends to every antipode of $\xi$.
\eprop

In the following we will study the branching behavior of horizontal geodesics in order to control the geometry of $\tits X$.

\blem\label{lem_trivialsplit}
Suppose that $\ga\subset\tits X$ is a piecewise $\delta$-horizontal geodesic joining two antipodes $\xi,\hat\xi\in\tits X$. If $\ga$ does not branch at an interior point,
then $\tits X$ splits off $\{\xi,\hat\xi\}$, $\tits X\cong\{\xi,\hat\xi\}\circ\pol\{\xi,\hat\xi\}$.
\elem

\proof
The assumption implies, via Proposition~\ref{prop_delnearalnear}, that $\hat\xi$ is the only antipode of $\xi$.
Since $X$ is geodesically complete, it has to agree with the parallel set $P(\xi,\hat\xi)$. This yields the claim.
\qed

\medskip

Our next goal is to show that the submetry $\delta:\tits X\to\Delta$ is transnormal (Proposition~\ref{prop_transnormal}).

To simplify notation we will make use of the following. For points $\xi\in\tits X$ and $o\in X$ we will denote the
corresponding direction $\log_o(\xi)$ by $\xi_o$.

\medskip

For a point $\xi\in\tits X$ denote by $N_\xi\subset\tits X$ the points which are $\delta$-near to $\xi$.
The subset $H_\xi\subset \Si_\xi\tits X$ denotes the {\em $\delta$-horizontal} directions, i.e. those directions where a $\delta$-horizontal geodesic starts.

Similarly, for a point $o\in X$ we denote by $N^o_v\subset\Si_oX$ the points which are $\delta_o$-near to $v$, and by $H_o\subset\Si_o X$
the $\delta_o$-horizontal directions. If $f:\tits X\to Z$ is a surjective 1-Lipschitz map to a geodesically complete CAT(1) space  $Z$,
the we we denote by $\delta^f: Z\to\Delta$ the induced submetry provided by Corollary~\ref{cor_imagedistsets}. 
If $f(x)=y$, then $H^f_y\subset\Si_y Z$ denotes the space of $\delta^f$-horizontal directions. 

If $Z$ is a CAT(1) space with a point $z\in Z$, then we call a subspace $H\subset\Si_z Z$ {\em almost symmetric},
if for every point $v\in H$ all genuine antipodes $\hat v\in Z$ are again contained in $H$.

\blem\label{lem_alsym}
Let $f:\tits X\to Z$ be a surjective 1-Lipschitz map  to a geodesically complete CAT(1) space $Z$ with $f(\xi)=x$.
If $H^f_{x}\subset\Si_{x} Z$ is almost symmetric, then $H_\xi$ is almost symmetric as well.
\elem

\proof
Let $\eta\in\tits X$ be $\delta$-near to $\xi$ and set $y=f(\eta)$.
Extend the geodesic $yx$ to an antipode $\hat y\in Z$ of $y$.
Denote by $\hat\eta\in \tits X$ a inverse image of $\hat Z$.
Since $f$ is 1-Lipschitz, Proposition~\ref{prop_delnearalnear} implies $\pi=|\hat\eta,\xi|+|\xi,\eta|$.
Thus, $f$ maps $\hat\eta\xi$ isometrically onto $\hat y x$. Since $H^f_x$ is almost symmetric,
there is a point $z\in x\hat y$ which is $\delta^f$-near to $x$. Its inverse image $\zeta$ on $\xi\hat\eta$
is then $\delta$-near to $\xi$. Since $\tits X$ is horizontally geodesically complete (Proposition~\ref{prop_delnearalnear}), it follows that
any geodesic extension of $\eta\xi$ is initially $\delta$-horizontal. Therefore $H_\xi$ is almost symmetric. 
\qed

\blem[{\cite[Lemma~8.2]{Ly_rigidity}}]
For $\xi\in\tits X$, $\hat \xi\in\Ant(\xi)$ let $\ga$ be a geodesic from $\xi$ to $\hat \xi$
with in- and outcoming directions $v$ and $w$, respectively. If $v\in H_\xi$ and $\hat w$
is a genuine antipode of $w$, then $\hat w\in H_{\hat \xi}$.
\elem

\blem[{\cite[Lemma~8.3]{Ly_rigidity}}]
If $H_\xi$ is almost symmetric, $\hat \xi\in \Ant(\xi)$, then $H_{\hat \xi}$ is almost symmetric too and there is a natural
isometry $I:H_\xi\to H_{\hat \xi}$.
\elem

\bcor[{\cite[Corollary~8.4]{Ly_rigidity}}]\label{cor_alsym}
If for $\xi\in\tits X$ the space of horizontal directions $H_\xi$ is almost symmetric, then for each 
$\zeta\in A_\xi=\bigcup_{j\in\N}\Ant^j(\xi)$ the subset $H_\zeta\in\Si_\zeta\tits X$ is almost symmetric as well.
Moreover, $H_\xi$ and $H_\zeta$ are isometric.
\ecor

\bprop\label{prop_transnormal}
Let $X$ be a locally compact and geodesically complete CAT(0) space, then the submetry $\delta:\tits X\to\Delta$
is transnormal. More precisely, for every pair of
$\delta$-near points $\xi$ and $\eta$ in $\tits X$, any extension of the segment $\xi\eta$
is piecewise $\delta$-horizontal.  
\eprop

\proof
For every point $o\in X$ the space of directions $\Si_o X$ is a compact geodesically complete CAT(1)
space. Thus, $\Si_oX$ admits many surjective 1-Lipschitz maps to round spheres \cite[Lemma~2.2]{Ly_rigidity}.
Via logarithmic maps we obtain surjective 1-Lipschitz maps from $\tits X$ to round spheres.
Let $f:\tits X\to S^n$ be such a map. For $\xi\in\tits X$ the set $f(A_\xi)$ is a dense subset of the fiber $f(B_\xi)$
of the induced submetry $\delta^f:S^n\to\Delta$. By \cite[Lemma~14.2]{Ly_rigidity}, we find a point $z=f(\zeta)\in f(A_\xi)$, such that the 
horizontal space $H^f_z$ is a round sphere. From Lemma~\ref{lem_alsym} and Corollary~\ref{cor_alsym} we conclude that $H_\xi$
is almost symmetric. This implies the claim. 
\qed

%

\blem\label{lem_uniformbranch}
Let $\xi\in\tits X$ be a point. For any $\delta$-horizontal direction $\hat v\in H_\xi$ exist times $0=t_0<t_1<\ldots<t_k=\pi$ 
with the following property.
For any genuine antipode $v$ of $\hat v$ and  any geodesic $\ga$ of length $\pi$ in the direction $v$, the point 
$\ga(t_i)$ is $\delta$-near to $\ga(t_{i+1})$.
\elem

\proof
By Lemma~\ref{lem_2flat}, there exists a geodesic $\hat\ga$ of length $\pi$ in the direction $\hat v$.
By Proposition~\ref{prop_transnormal}, there exist times $0=t_0<t_1<\ldots<t_k=\pi$  such that $\hat\ga(t_i)$ is $\delta$-near to $\hat\ga(t_{i+1})$.
We see that $\ga(t)\in B_{\hat\ga(t)}$. Thus if $\hat\ga(t)$ is $\delta$-near to $\hat\ga(s)$, then $\ga(t)$ is $\delta$-near to $\ga(s)$.
\qed

\blem\label{lem_unihor}
For $\xi\in\tits X$ there exists $\eps>0$ such that for each $v\in H_\xi$, there is a unique geodesic of length $\eps$ in the direction $v$, and this geodesic
 is $\delta$-horizontal.
\elem

\proof
Choose a complete geodesic $c$ in $X$ asymptotic to $\xi$ and let $o$ be a point on $c$.
By \cite[Proposition~5.3]{Ly_rigidity}, there exists a positive $\eps=\eps(\delta(\xi))$ such that every $\delta_o$-horizontal geodesic in  $\Si_oX$
starting in $\xi_o$ lies on a $\delta_o$-horizontal geodesic of length $\eps$.

By definition of $H_\xi$, there is a $\delta$-horizontal geodesic in $\tits X$ in direction $v$. By Lemma~\ref{lem_2flat}, we find a flat half-plane $H\subset X$ whose boundary is $x$ and whose ideal boundary is a 
geodesic $\ga$ in the direction $v$. Then $\log_o$ maps $\ga$ isometrically to a geodesic $\ga_o$ in $\Si_oX$.
By Corollary~\ref{cor_fact}, $\ga_o$ is initially $\delta_o$-horizontal and therefore $\delta_o$-horizontal up to time $\eps$. 
Since $\delta$ factors as $\delta_o\circ\log_o$, we see that $\ga$ is $\delta$-horizontal up to time $\eps$.
By Proposition~\ref{prop_sph}, any geodesic in direction $v$ has to agree with $\ga$ up to time $\eps$.
\qed

\blem\label{lem_nearclosed}
For every point $\xi\in\tits X$ the set $N_\xi\subset\geo X$ is closed with respect to the cone topology.
\elem

\proof
Let $\eta_k\to\eta$ be a convergent sequence  in $N_\xi$. By Lemma~\ref{lem_unihor}, we may assume $|\xi,\eta_k|\geq\eps(\xi)>0$.
Choose a complete geodesic $c$ in $X$ asymptotic to $\xi$ and let $o$ be a point on $c$.
 By Lemma~\ref{lem_2flat},
there are flat half-planes $H_k$ in $X$ with $\D H_k=c$ and $\eta_k\in\geo H_k$.
After passing to a subsequence, we obtain a limit flat half-plane $H$ with $\D H=c$ and $\eta\in\geo H$.
In particular, $|\xi,\eta|=\angle_o(\xi,\eta)$. Since $\delta_o$ is continuous, $\eta_o$ is $\delta_o$-near to $\xi_o$. 
The claim follows since $\delta$ factorizes as $\delta_o\circ\log_o$.
\qed

\blem\label{lem_nearconvex}
For every point $\xi\in\tits X$ the set $N_\xi\subset\tits X$ is convex.
\elem

\proof
Choose points $\eta, \eta'\in N_\xi$ at distance less than $\pi$ and denote by $\zeta$ their midpoint. 
It is enough to show that $\zeta$ is $\delta$-near to $\xi$.
Choose a complete geodesic $c$ in $X$ asymptotic to $\xi$ and let $o$ be a point on $c$.
Let $\rho$ be a parametrization of the geodesic ray $o\eta$. Set $\delta_t:=\delta_{\rho(t)}$, $\xi_t:=\xi_{\rho(t)}$
and $\zeta_t:=\zeta_{\rho(t)}$. By the Lune Lemma, there is a geodesic in $\tits X$ through $\zeta$ which connects the ideal 
boundary points of $c$.  In particular, there exists a flat half-plane $H_t$ in $X$ whose boundary $\D H_t$ is parallel to $c$, contains the point $\rho(t)$
and with $\zeta\in\geo H_t$. We conclude that $\angle_{\rho(t)}(\xi_t,\zeta_t)=|\xi,\zeta|$ holds for all $t\geq 0$.

By \cite[Lemma~8.1]{Ly_rigidity}, the set $N_{\xi_t}\subset\Si_{\rho(t)}X$
is convex. Hence, the midpoint $m_t\in \Si_{\rho(t)}X$ of $\eta_t$ and $\eta'_t$ lies in $N_{\xi_t}$.
For every $\eps>0$ we find $t_\eps>0$ such that $|\eta,\eta'|\leq\angle_{\rho(t)}(\eta,\eta')+\eps$ for $t\geq t_\eps$.
Thus,
\[\angle_{\rho(t)}(\eta,\zeta)+\angle_{\rho(t)}(\zeta,\eta')\leq|\eta,\eta'|\leq\angle_{\rho(t)}(\eta,\eta')+\eps\]
holds for $t\geq t_\eps$. From triangle comparison we conclude $\angle_{\rho(t)}(m_t,\zeta_t)\to 0$. 
Hence,
\begin{align*}
|\delta(\xi),\delta(\zeta)|&=|\delta_t(\xi_t),\delta_t(\zeta_t)|\geq |\delta_t(\xi_t),\delta_t(m_t)|-|\delta_t(m_t),\delta_t(\zeta_t)|\\
&\geq\angle_{\rho(t)}(\xi_t,m_t)-\angle_{\rho(t)}(m_t,\zeta_t)\\
&\geq\underbrace{\angle_{\rho(t)}(\xi_t,\zeta_t)}_{=|\xi,\zeta|}-2\underbrace{\angle_{\rho(t)}(m_t,\zeta_t)}_{\to 0}\to|\xi,\eta|.
\end{align*}
It follows that $\xi$ and $\zeta$ are $\delta$-near as required.
\qed

\subsection{Spherical splittings}\label{sec_splittings}

In this section we collect some splitting results which can be read off the branching behavior of $\delta$-horizontal
geodesics. First, let us recall the following general result.

\bthm[{\cite[Theorem~2.9]{BMS_affine}}]\label{thm_split}
Let $X$ be a locally compact CAT(0) space with finite-dimensional Tits boundary and where $\geo X$ contains a dense family of
top-dimensional round spheres. Suppose
$\geo X$ contains a non-empty proper closed convex subset $Z$ which is symmetric. 
Then $\tits X$ splits as $\tits X\cong Z\circ Z^\perp$.
\ethm

Note, by Corollary~\ref{cor_morsecyclesdense}, this applies to all locally compact CAT(0) spaces with geometric group actions. 
Thus in order to find spherical splittings of such spaces, one is on the lookout for closed convex symmetric sets at infinity.

\bprop\label{prop_split} 
Let $\hat\si=\si_1\circ\si_2$ be a splitting into round spheres.
Suppose that every pair of points $\xi_i\in\si_i$, $i=1,2$, the geodesic $\xi_1\xi_2$ is piecewise horizontal and does not branch
in $\tits X$. Then $\Ant(\si_i)$, the set of antipodes of $\si_i$, is a closed convex symmetric proper subset of $\geo X$. 
\eprop

\proof
Note first that if a point $\zeta\in\tits X$ has constant distance $\frac{\pi}{2}$ from every point in $\si_1$, then it has to lie in $\Ant(\si_2)$.
Indeed, by Corollary~\ref{cor_morseantipodes}, $\zeta$ has an antipode $\hat\zeta$ in $\hat\si$. By assumption, $\hat\zeta$ lies on a non-branching piecewise horizontal geodesic $\xi_2\xi_1$
with $\xi_i\in\si_i$. Then the segment $\hat\zeta\xi_1$ is initially horizontal and since $\xi_2\xi_1$ does not branch, we can extend 
$\hat\zeta\xi_1$ to $\zeta$. Thus $\hat\zeta=\xi_2\in\si_2$.

Now we show that $\Ant(\si_2)$ is closed. Let $\hat\zeta_k\to\hat\zeta$ be a convergent sequence in $\geo X$ with $\hat\zeta_k\in\Ant(\zeta_k)$, $\zeta_k\in\si_2$.
Let $\xi_1\in\si_1$ be an arbitrary point and denote by $\hat\xi_1\in\si_1$ its antipode.
Arguing as before, we find geodesics $\ga_k$ and $\hat\ga_k$ from $\zeta_k$ to $\hat\zeta_k$ through the point $\xi_1$ and $\hat\xi_1$, respectively.
Passing to limits, we obtain geodesics $\ga$ and $\hat\ga$ from a point $\zeta\in\si_2$ to $\hat\zeta$ through the point $\xi_1$ and $\hat\xi_1$, respectively.
By semicontinuity of the Tits metric, we conclude $|\hat\zeta,\xi_1|=|\hat\zeta,\hat\xi_1|=\frac{\pi}{2}$. Therefore, $\hat\zeta\in\Ant(\zeta)\subset \Ant(\si_2)$. 

Before we show symmetry, observe that if $\hat\zeta$ is antipodal to a point $\zeta\in\si_2$,
and $\xi_1\in\si_1$ is arbitrary, then the geodesic $\hat\zeta\xi_1$ is piecewise horizontal and does not branch.
Indeed, it is contained in a round $1$-sphere which also contains the points $\zeta$ and $\hat\xi_1$.
Since $\zeta\hat\xi_1$ does not branch, neither does $\hat\zeta\xi_1$. 
Thus, if $\hat\zeta'$ is another point in $\Ant(\si_2)$ which is antipodal to $\hat\zeta$, then there is a geodesic 
from $\hat\zeta$ to $\hat\zeta'$ through the point $\xi_1$. In particular, $|\hat\zeta',\xi_1|=\frac{\pi}{2}$.
Since $\xi_1\in\si_1$ was arbitrary, we conclude $\hat\zeta'\in\Ant(\si_2)$ as required.

At last, let us show convexity. Let $\hat\zeta_2$ and $\hat\zeta_2'$ be points at distance less than $\pi$ and antipodal to points $\zeta_2$ and $\zeta_2'$ in $\si_2$,
respectively. For every point $\xi_1\in\si_1$ the Lune Lemma provides a spherical lune spanned by $\xi_1,\hat\xi_1$
 and $\hat\zeta_2,\hat\zeta_2'$. In particular, the segment $\hat\zeta_2\hat\zeta_2'$  has constant distance $\frac{\pi}{2}$ from $\xi_1$.
Again, since $\xi_1\in\si_1$ was arbitrary, we conclude $\hat\zeta_2\hat\zeta_2'\subset\Ant(\si_2)$ as required.
\qed
\medskip

The following auxiliary result, which guarantees a particular branching pattern along horizontal geodesics, will be used later in order to find
join decompositions as required by Proposition~\ref{prop_split}.

\blem\label{lem_split}
Suppose that there exists a point $\xi\in\tits X$ and a direction $v\in H_\xi$ such that $\Ant(v)$ contains a round 1-sphere $\si$ of genuine directions.
Then any geodesic in direction $v$ can only contain branch points at times which are multiples of $\frac{\pi}{2}$.
\elem

\proof
Choose an antipode $\hat\xi$ of $\xi$. By Proposition~\ref{prop_transnormal}, the circle $\si\subset\Si_\xi\tits X$ consists of $\delta$-horizontal directions. 
By the Lune Lemma, the union of geodesics from $\xi$ to $\hat\xi$ in the directions of $\si$
form a round 2-sphere $\bar\si\subset\tits X$. By Lemma~\ref{lem_uniformbranch}, there exist 
times $0=t_0<t_1<\ldots<t_k=\pi$ such that any geodesic $\ga$ connecting $\xi$ to $\hat\xi$ in a  direction of $\si$ has the property that 
$\ga(t_i)$ is $\delta$-near to $\ga(t_{i+1})$. Let $\ga$ be such a geodesic. 
The claim follows if we can show that
 $\ga$ does not branch for $t\neq\frac{\pi}{2}$.
More precisely, we need to rule out that $\ga$ or its inverse parametrization branch at a time $t_{i_0}\neq\frac{\pi}{2}$. 
Since the setting is symmetric, we can assume that $\ga$ itself branches. Moreover, we may assume $t_{i_0}>\frac{\pi}{2}$.
Indeed, if $\hat\ga$ denotes the unique geodesic in $\hat\si$, starting in $\xi$ and such that $\ga\cup\hat\ga$ form a closed geodesic,
then $\hat\ga$ branches at $\hat t=\pi-t_{i_0}$. So if $t_{i_0}<\frac{\pi}{2}$, then we consider $\hat\ga$ instead.
Now let $\hat\xi'$ denote the antipode of $\xi$ obtained by prolonging the segment $\xi\ga(t_{i_0})$ along the branch.
Denote by $w\in\Si_\xi\hat\si$ the starting direction of $\ga$. Let $w'\in\Si_\xi\hat\si$ be a direction close to $w$.
We show that the piecewise $\delta$-horizontal geodesic $\ga'$ from $\xi$ to $\hat\xi'$ in the direction $w'$ has to branch in $t_{i_0}$.

\begin{center}
\includegraphics[scale=0.3,trim={4.5cm 1.5cm 0cm 7cm},clip]{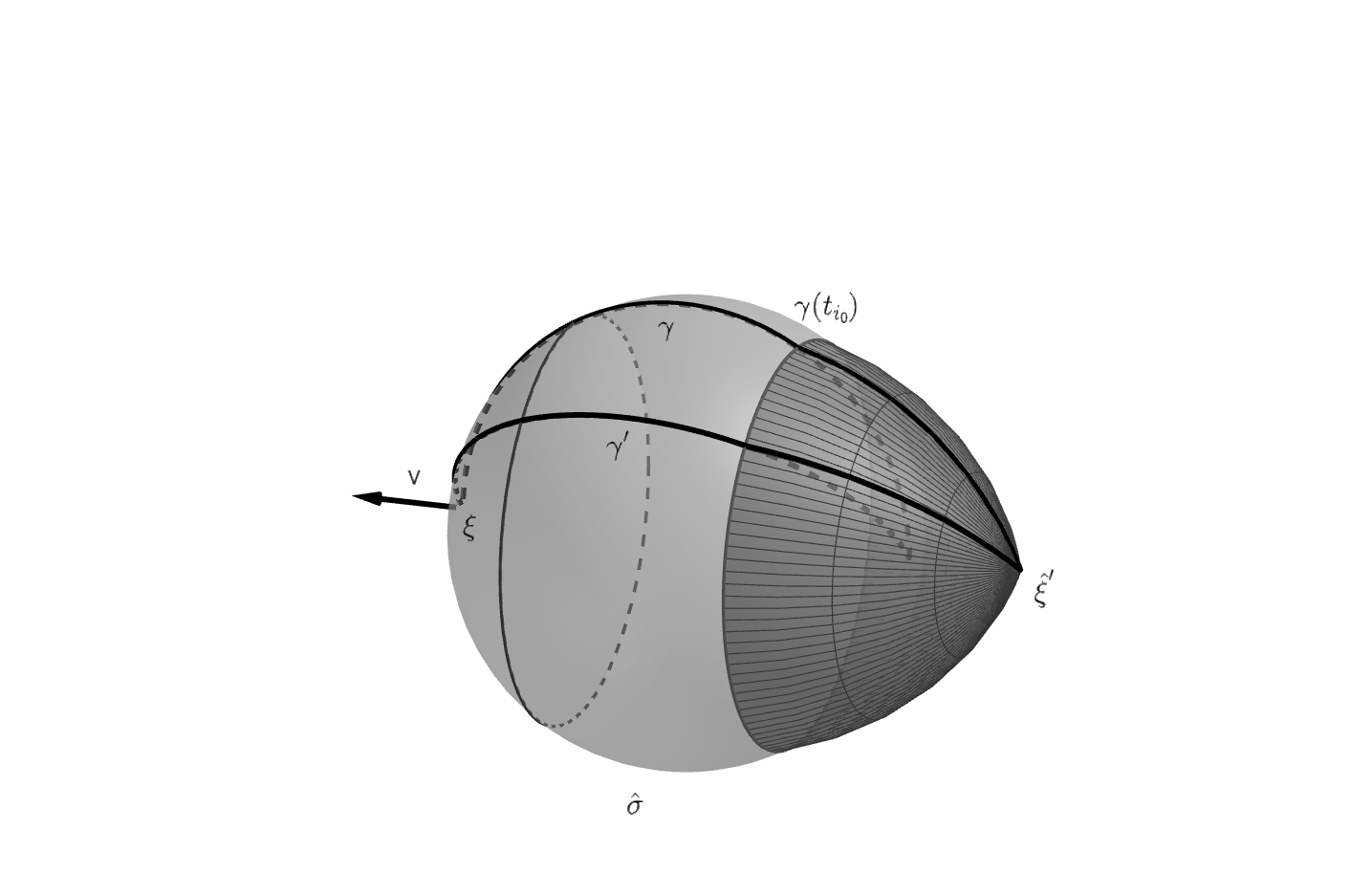}
\end{center}

By the Lune Lemma, we find two spherical lunes $\la$ and $\la'$ spanned by the directions $w$ and $w'$ and corresponding to the points $\hat\xi$ and $\hat\xi'$. If $\ga'$ branches at time $t_j$, then the segment $\ga(t_{i_0})\ga'(t_j)$ is contained in $\la\cap\la'$. Thus the directions between $w$ and $w'$
realize every branch time between $t_{i_0}$ and $t_j$. Since branching is only possible at the finitely many times $t_i$, we conclude $t_j=t_{i_0}$.
By continuity, we infer that every geodesic from $\xi$ to $\hat\xi'$ in a direction of $\si$ has to branch at $t_{i_0}$.
Let $\hat\ga$ denote the particular geodesic which starts in direction $\hat w\in\Si_\xi\hat\si$ antipodal to $w$.
Then the points  $\ga(t_{i_0})$ and $\hat\ga(t_{i_0})$ are joined by two geodesics, one through $\hat\xi$ and one through $\hat\xi'$.
This is a contradiction since $|\ga(t_{i_0}),\hat\ga(t_{i_0})|=2(\pi-t_{i_0})<\pi$.
\qed

\section{Morse flats and higher rank}\label{sec_final}

\subsection{Tits boundary of a periodic Morse flat}

Throughout this section,  we assume that $X$ is a locally compact  CAT(0) space which contains a periodic Morse $(n+1)$-flat $\hat F$ with $\hat\si:=\tits \hat F$.
Recall from \cite{St_rr} that
a round sphere $\si\subset\hat\si$ is called {\em singular}, if there exists a round hemisphere $\tau^+\subset\tits X$
with $\D\tau^+=\si$ and  $\tau^+\cap\hat\si=\si$ such that the union $\tau^+\cup\tau^-$ is a round sphere for every 
round hemisphere $\tau^-\subset\hat\si$ with $\D\tau^-=\si$. 

\blem\label{lem_goodhemi}
Let $\tau\subset\tits X$ be a round $n$-hemisphere with $\tau\cap\hat\si=\D\tau$. Then 
$\tau\cup\hat\si$
is a convex subset in $\tits X$ which splits off $\D\tau$. Thus $\D\tau$
is a singular hypersphere in $\hat\si$.
\elem

\proof
Let $\zeta$ be the center of $\tau$.
Denote by $\zeta^\pm\in\hat\si$ the two points perpendicular to $\D\tau$.
Since $\hat\si$ is the Tits boundary of a Morse flat, $\zeta$ has to be antipodal to both points $\zeta^\pm$ (Proposition~\ref{prop_morseinj}).
Thus, the claim follows from the Lune Lemma.
\qed

\bthm[{\cite[Theorem~4.7]{St_rr}}]\label{thm_singfinite}
Let $X$ be a locally compact CAT(0) space.  Suppose that $X$ contains a periodic Morse $(n+1)$-flat $\hat F$ with $\hat\si:=\tits \hat F$. Then 
the complement in $\tits F$ of the set of regular points $\hat\si_{reg}$ can be covered by a finite set of singular spheres of positive codimension.
\ethm


\medskip

%

In the following we denote by $\si^1,\ldots,\si^s\subset\hat\si$ the finitely many singular hyperspheres provided by 
Theorem~\ref{thm_singfinite}.
For a point $\xi\in\hat\si_{reg}$ we define 
\[C_\xi:=\bigcap_{\xi\in\tilde\si}\tilde\si\subset\hat\si\] 
as the intersection of all round $n$-spheres $\tilde\si$ in $\tits X$ containing $\xi$.
By Theorem~\ref{thm_singfinite}, $C_\xi$ contains a neighborhood of $\xi$ and therefore has dimension $n$.
Our next aim is to show that $C_\xi$ agrees with the closure of one of the components of $\hat\si\setminus\bigcup_{i=1}^s\si^i$, in
particular it's an $n$-dimensional spherical polytope. 
We need some preparation.

Recall that by Rademacher's theorem, almost every point $p$ in the boundary $\D C$ of an $n$-dimensional closed convex subset $C$ in $S^n$ 
has a Euclidean tangent space, $T_p\D C\cong\R^{n-1}$. 

\bprop\label{prop_bomb}
Let $X$ be a locally compact CAT(0) space and $\hat F\subset X$ a periodic Morse $(n+1)$-flat with $\hat\si:=\tits \hat F$.
Let $\hat F'\subset X$ be a second $(n+1)$-flat with $\hat\si':=\tits \hat F'$.
Suppose that $\al=\hat\si\cap\hat\si'$ is non-empty and of dimension $n$.
Then for every point $\xi\in\D\al$ with $T_\xi\D\al\cong \R^{n-1}$, there exists a singular hypersphere $\si\subset\hat\si$ 
which 
is tangent to $\al$ at $\xi$. 
More precisely,  $\xi\in\si$ and $\al\cap\si\subset\D\al$.
\eprop

\proof
Since $X$ is locally compact and $\hat F$ is Morse, there exists $R>0$, such that $\hat F$
does not bound a flat strip of width $R$. Set $C_R:=\{x\in\hat F:\ |x,\hat F'|\leq R\}$
and pick a base point $o\in \D C_R$. Note that $C_R$ is a closed convex subset of $\hat F$.
If $C_o(\al)$ denotes the cone from $o$ over $\al$, then $C_o(\al)\subset C_R$ and $C_R$
has sublinear distance from $C_o(\al)$ in  the sense 
\[\frac{|C_o(\al)\cap B_r(o),C_R\cap B_r(o)|}{r}\to 0.\]
Thus, for $\xi\in\D\al$, we can find a sequence $x_k\in\D C_R$ such that $x_k\to\xi$. Denote by $x_k'\in\hat F'$
the closest point to $x_k$. Thus $|x_k,x_k'|=R$. Now we use elements $\ga_k$ in the stabilizer of $\hat F$ to move $x_k$
to a fixed compact subset of $\hat F$. We may assume that each  $\ga_k$ acts as a translation on $\hat F$ and therefore fixes $\hat\si$
pointwise. Passing to a subsequence, we obtain limits $\ga_k(x_k)\to x_\infty$ and $\ga_k(\hat F',x_k')\to(\hat F_\infty', x_\infty')$. 
By construction, $|x_\infty,\hat F'_\infty|=|x_\infty,x'_\infty|=R$. In particular, $\hat F'_\infty$ cannot be parallel to $\hat F$.
We claim that $\hat F'_\infty$ is asymptotic to $\hat F$ along an $(n+1)$-dimensional flat half-space.
To see this, note first that $\angle_o(x_k,\xi)\to 0$ and therefore $\angle_{x_k}(o,\xi)=\pi-\angle_o(x_k,\xi)\to \pi$.
Thus we have the pointed convergences $\ga_k(x_k\xi,x_k)\to(x_\infty\xi,x_\infty)$
and $\ga_k(x_k o,x_k)\to(x_\infty\hat\xi,x_\infty)$ where $\hat\xi$ is the antipode of $\xi$ in $\hat\si$.
We denote by $c$ the complete geodesic through $x_\infty$ which is asymptotic to the points $\xi,\hat\xi$.
Since the ray $x_k\xi$ and the segment $x_k o$ are contained in $C_R$, we conclude $c\subset \bar N_R(\hat F_\infty')$.
(Thus $c$ is parallel to $\hat F_\infty'$ and since $|x_\infty,\hat F'_\infty|=R$, it has distance $R$ from $\hat F_\infty'$.)
In particular, $\xi,\hat\xi\in\tits\hat F_\infty'$. Since $\al$ is fixed by $\ga_k$ and contained in $\hat\si'$, we also have $\al\subset\tits\hat F_\infty'$.
Therefore, $\tits\hat F_\infty'$ has to contain the convex hull of $\hat\xi$ and $\al$ in $\hat\si$.
By assumption, the tangent space to $\D\al$ at $\xi$ is $\R^{n-1}$ (and therefore the tangent space to $\al$
at $\xi$ is an $n$-dimensional Euclidean half-space). It follows that this convex hull has to contain the round $n$-hemisphere 
$\hat\tau\subset\hat\si$ determined by $\hat\xi$ and $\al$. Since $\hat F'_\infty$ is not parallel to $\hat F$, we conclude 
$\tits\hat F_\infty'\cap\hat\si=\hat\tau$. Thus $\si=\D\hat\tau$ is a singular hypersphere as claimed.
\qed

\bcor\label{cor_poly}
Let $X$ be a locally compact CAT(0) space and $\hat F\subset X$ a periodic Morse $(n+1)$-flat with $\hat\si:=\tits \hat F$.
Let $\hat F'\subset X$ be a second $(n+1)$-flat with $\hat\si':=\tits \hat F'$.
Suppose that $\al=\hat\si\cap\hat\si'$ is non-empty and of dimension $n$. Then
there exists finitely many singular hyperspheres $\{\si^1,\ldots,\si^k\}$ in $\hat\si$  such that
$\al$ is the closure of a component of
$\hat\si\setminus\bigcup_{i=1}^k\si^i$. In particular,
 $\al$ is an $n$-dimensional spherical polytope.
\ecor

\proof
By Proposition~\ref{prop_bomb}, almost every point in $\D\al$ is tangent to one of the singular hyperspheres $\si^i$. 
Let $\{\si^1,\ldots,\si^k\}$ denote those which are tangent to some point in $\D\al$. 
 It follows that $\al$ has to agree with the closure of a component of $\hat\si\setminus\bigcup_{i=1}^k\si^i$.
\qed

\bcor\label{cor_regpoly}
There exists finitely many round singular hyperspheres $\{\si^1,\ldots,\si^s\}$ in $\hat\si$ such that
for every point 
$\xi\in\hat\si_{reg}$, the set
$C_\xi$ is the closure of one of the components of
$\hat\si\setminus\bigcup_{i=1}^s\si^i$. In particular,
there are only finitely many distinct $C_\xi$ and
each $C_\xi$ is an $n$-dimensional spherical polytope.
\ecor

\proof
Choose a point
$\xi\in\hat\si_{reg}$.
Recall that $C_\xi$ is the intersection of all round $n$-spheres in $\tits X$ which contain the point $\xi$.
By Corollary~\ref{cor_poly}, every finite such intersection has the claimed structure.
Since by Theorem~\ref{thm_singfinite} there are only finitely many singular hyperspheres, 
the set $C_\xi$ can be obtained as a finite intersection of round $n$-spheres.
At last, observe that if $C_\xi$ is bounded by the hyperspheres $\si^1,\ldots,\si^j$, then none of the other singular hyperspheres $\si^i$
can cut through $C_\xi$, since, by Lemma~\ref{lem_goodhemi}, a corresponding hemisphere $\tau_i$ spans a round $n$-sphere with each of the components of $\hat\si\setminus\si^i$. 
The claim follows.
\qed

%

\subsection{Periodic Morse flats and symmetric sets}

In this section we consider a locally compact CAT(0) space $X$ with a geometric group action $\Ga\acts X$
and such that every geodesic in $X$ lies in an $(n+1)$-flat, $n\geq 1$. We
 will show that a 
$\Ga$-periodic Morse $(n+1)$-flat $\hat F\subset X$ leads to a non-empty proper closed symmetric $\Ga$-invariant subset in $\geo X$ and therefore, by Proposition~\ref{prop_submetry}, to a non-trivial submetry of the
Tits boundary. The proof is dynamical and we will freely use the notation and results from Section~\ref{subsec_dyn}.

 As before, we will set $\hat\si=\tits \hat F$ and denote the finitely many singular hyperspheres in $\hat\si$ by $\si^1,\ldots,\si^s$.
Moreover, we define a {\em chamber (in $\hat\si$)} to be  the closure of one of the components of $\hat\si\setminus\bigcup_{i=1}^s\si^i$.
By Corollary~\ref{cor_regpoly}, for each regular point $\xi\in\hat\si_{reg}$ the intersection of all round $n$-spheres in $\tits X$ which contain $\xi$ is a chamber.

\blem\label{lem_antipodesinflats}
Let $C\subset\hat\si$ be a chamber and let $\xi$ be a point in the interior of $C$. If $\hat\xi\in\tits X$
is antipodal to $\xi$,
then there exists an $(n+1)$-flat $F\subset X$ such $C$ and $\hat\xi$ are contained in $\tits F$.
\elem

\proof
If $\hat\xi$ is visually antipodal to $\xi$, then the claim follows since every geodesic in $X$
lies in an $(n+1)$-flat and $C$ is the intersection of all round $n$-spheres in $\tits X$
which contain a regular point from $C$ (Corollary~\ref{cor_regpoly}).
By Lemma~\ref{lem_dense}, there exists a sequence $(\hat\xi_k)$ of visual antipodes of $\xi$ with $\hat\xi_k\to\hat\xi$.
Denote by $F_k$ a $(n+1)$-flat such that $\si_k:=\tits F_k$ contains $C$ and $\hat\xi_k$.
Let $f_k:S^n\to\tits X$ be an isometric embedding with image $\si_k$.
Passing to a subsequence, we can assume that $f_k$ converges pointwise with respect to the cone topology to a 1-Lipschitz map
$f:S^n\to\tits X$. Since the image of $f$ contains $C$ and $\hat\xi$, it has to be a round $n$-sphere $\si$.
Because $C\subset\hat\si$, we see that $\si$ cannot bound a round hemisphere in $\tits X$.
Thus \cite[Proposition~2.1]{Leeb} implies the claim.
\qed
\medskip

In order to show that chambers are incompressible we will make use of the following.

\blem\label{lem_regularantipodes}
Let $C\subset\hat\si$ be a chamber and $\ga\in\Ga$. Suppose that $\xi\in\ga C$
is a point such that $\ga^{-1}\xi$ is regular in $C$. Let $\hat\xi\in\hat\si$ be an antipode of $\xi$.
Then $\hat\xi$ is regular as well. More precisely, if $s>0$ is such that $B_s(\xi)\subset \ga C$,
then there exists a chamber $\hat C\subset\hat\si$ such that $B_s(\hat\xi)\subset \hat C$.
\elem

\proof
Choose $s>0$ as in the assumption and
let $\hat\xi'\in B_s(\hat\xi)$. Extend the geodesic $\hat\xi'\hat\xi$ inside $\hat\si$ to a point $\hat\xi''\in\hat\si$
such that $\hat\xi\hat\xi''$ lies in a single chamber $\hat C\subset\hat\si$. Choose a regular point $\hat\eta\in\hat C\cap B_s(\hat\xi)$.
Extend the geodesic $\hat\eta\xi$ to antipode $\eta$. Since $|\eta,\xi|>\pi-s$ and $B_s(\xi)\subset \ga C$, we have $\eta\in \ga C$.
By Lemma~\ref{lem_antipodesinflats}, we find a round $n$-sphere $\si\subset\tits X$ which contains the chamber $\hat C$ and the point $\eta$, and therefore
also the translated chamber $\ga C$. In particular, we see that $\xi$ is the only antipode of $\hat\xi$ in $\ga C$ and $\si$ is given by the union of all geodesics
from $\xi$ to $\hat\xi$. Now we extend the geodesic $\hat\xi'\hat\xi''$ inside $\si$ to a geodesic $c'$ of length $\pi$ ending in an antipode $\xi'$ of $\hat\xi'$.
Since $|\hat\xi',\hat\xi|<s$ we see that $\xi'\in B_s(\xi)\subset\ga C$. By Lemma~\ref{lem_antipodesinflats}, we find a round $n$-sphere $\si'$
in $\tits X$ which contains the translated chamber $\ga C$ and the point $\hat\xi'$. Again,  $\xi'$ is the only antipode of $\hat\xi'$ in $\ga C$
and therefore $\si'$ is the union of all geodesics from $\xi'$ to $\hat\xi'$. In particular, $\si'$ contains the geodesic $c'$ and therefore also the point $\hat\xi$.
But now both round $n$-spheres $\si$ and $\si'$ contain the point $\hat\xi$ and the translated chamber $\ga C$. Since $\xi$ lies in the interior of $\ga C$
and $\tits X$ is CAT(1), this implies $\si=\si'$. Consequently, $B_s(\hat\xi)\subset\si$. Since $\si$ is a round $n$-sphere, we must have 
$B_s(\hat\xi)\subset\hat\si$ as claimed.
\qed

\blem\label{lem_chamincomp}
Let $C\subset\hat\si$ be a chamber, then $C$ is incompressible.
\elem

\proof
Let $\eta$ be a regular point in $C$ and assume $(\ga_k)$ is a sequence in $\Ga$
with $\eta_k\to\xi$ with $\eta_k=\ga_k\eta$. We need to show that a limit of $\ga_k C$ with respect to the cone topology is isometric to $C$. 

Let $\hat\xi$ be an antipode of $\xi$ in $\hat\si$ (Corollary~\ref{cor_morseantipodes}).
Since the diameter of $\tits X$ is equal to $\pi$, semi-continuity of the Tits metric implies  $|\hat\xi,\eta_k|\to\pi$.
Extend the geodesic $\eta_k\hat\xi$ inside $\hat\si$ to an antipode $\hat\xi_k\in\hat\si$. Then $\hat\xi_k\to\hat\xi$.
By Lemma~\ref{lem_regularantipodes}, if $s_0>0$ is such that $B_{s_0}(\eta)\subset C$, then $B_{s_0}(\hat\xi_k)\subset \hat\si$. 
Therefore $\hat\xi$ is regular and if $\hat C\subset\hat\si$ denotes the chamber containing $\hat\xi$, then by Lemma~\ref{lem_antipodesinflats},
we find a round $n$-sphere $\si_k\subset\tits X$ which contains $\hat C$ and $\ga_k C$. Let $f_k:S^n\to\tits X$ be an isometric embedding
with image $\si_k$. Passing to a subsequence, we may assume that $f_k$ converges pointwise with respect to the cone topology to a $1$-Lipschitz
map $f:S^n\to\tits X$. Since the image of $f$ has to contain $\hat C$ and the point $\xi$, we see that $f$ is an isometric embedding.
In particular, any limit with respect to the cone topology of $\ga_k C$ is isometric to $C$ as claimed.
\qed

%
%

\blem
The chambers of $\hat\si$ are precisely the maximal incompressible sets in $\hat\si$.
\elem

\proof
Let $\si^1$ be a singular hypersphere in $\hat\si$ and let $\tau\subset\tits X$ be a round $n$-hemisphere
with $\tau\cap\hat\si=\si^1$. Denote by $\tau^\pm\subset\hat\si$ the two round $n$-hemispheres
determined by $\si^1$. Since $\hat\si$ is the boundary of a Morse flat, the two round $n$-spheres $\si^\pm=\tau\cup\tau^\pm$
do not bound round hemispheres in $\tits X$. Thus by \cite[Proposition~2.1]{Leeb}, there are $(n+1)$-flats $F^\pm\subset X$
with $\tits F^\pm=\si^\pm$. Let $\zeta\in\tits X$ and $\zeta^\pm\in\tits X$ be the center of $\tau$ and $\tau^\pm$, respectively.
If $\om\in\beta\Ga\setminus\Ga$ pulls from $\zeta$, then, by Lemma~\ref{lem_pullflat}, $T^\om\si^+$ and $T^\om\si^-$ are  round $n$-spheres. Moreover,  
by $\pi$-convergence, $T^\om\zeta^+=T^\om\zeta^-$ and   the round hemispheres $\tau^+$ and $\tau^-$
get folded to the same round hemisphere in $T^\om\si^+=T^\om\si^-$.  Since the chambers in $\hat\si$ are bounded by the singular hyperspheres, this shows that no incompressible 
set can strictly contain a chamber. By Lemma~\ref{lem_chamincomp}, each chamber is incompressible and the proof is complete. 
\qed
\medskip

Following \cite{GuSw_trans}, we will denote the non-degenerated maximal incompressible sets in $\tits X$ by $\mathcal{I}(\Ga)$.
By Lemma~\ref{lem_morseretract}, the maximal dimension of an incompressible set in $\tits X$ is at most $n$. Thus we have

\blem\label{lem_strongfold}
Let $C_{\max}$ be a chamber of maximal volume in $\hat\si$. Then there exists $\nu\in\beta\Ga$
such that $T^\nu \tits X= T^\nu C_{\max}$.
\elem

\proof
This is part of the Strong Folding Lemma \cite[Theorem~D]{GuSw_trans}.
\qed

\bcor\label{cor_alliso}
All chambers of $\hat\si$ are isometric.
\ecor

\proof
Since every chamber is incompressible, the map $T^\nu$ from Lemma~\ref{lem_strongfold},
restricts to an isometric embedding $C\hookrightarrow T^\nu C_{\max}$ for every chamber $C\subset\hat\si$.
Thus the inverse image of $\D T^\nu C_{\max}$ in $\hat\si$ is precisely the union of the singular hyperspheres $\bigcup_{i=1}^s\si^i$.
Hence each restriction $C\hookrightarrow T^\nu C_{\max}$ is surjective.
\qed

\bcor\label{cor_chamberscover}
For every point $\xi\in\tits X$ there exists $\om\in\beta\Ga$ and a chamber $C\subset\hat\si$ such that $\xi\in T^\om C$. In particular, $\tits X$ is covered by maximal incompressible subsets isometric to a chamber of $\hat\si$.
\ecor

\proof
Since $\hat\si$ is covered by chambers, the claim follows because the orbit  $\Ga\hat\si$ is dense (Corollary~\ref{cor_morsecyclesdense}).
\qed

\blem\label{lem_center}
Every chamber in $\hat\si$ has an incenter unless $X$ is flat.
\elem

\proof
A chamber lies in a hemisphere of $\hat\si$ unless the singular set in $\hat\si$ is empty.
But then $\tits X=\hat\si$ and $X$ is isometric to $\R^{n+1}$.
\qed

\bprop\label{prop_invarantset}
If $X$ is not flat, then there exists a non-empty proper closed symmetric  $\Ga$-invariant subset $A\subset\geo X$. 
\eprop

\proof
Denote by $T^\nu$ the folding map $\tits X\to T^\nu C_{\max}$ from Lemma~\ref{lem_strongfold}.
Let $\mu$ be the incenter of $T^\nu C_{\max}$ (Lemma~\ref{lem_center}).
Let $A$ be the inverse image of $\mu$ under $T^\nu$, i.e. the incenters of maximal incompressible subsets of $\tits X$.
Then $A$ is a non-empty proper closed $\Ga$-invariant set. We claim that $A$ is also symmetric.
Let $\xi$ be in $A$ and choose an antipode $\hat\xi$ of $\xi$. 
By Corollary~\ref{cor_chamberscover}, there exist sequences $(\ga_k)$ and $(\hat\ga_k)$ in $\Ga$
and chambers $C$ and $\hat C$ in $\hat\si$ such that $\ga_k C\to C'$, $\hat\ga_k\hat C\to\hat C'$ and $\xi\in C'$, $\hat\xi\in\hat C'$.
Since $\xi\in A$, the incenters $\xi_k$ of $\ga_k C$ converge to $\xi$. To prove the claim it is enough to show that the incenters $\hat\xi_k$
of $\hat\ga_k\hat C$ converge to $\hat\xi$. Let $\zeta_k\in\hat\ga_k \hat C$ be regular and such that $\zeta_k\to\hat\xi$.
Then $|\xi_k,\zeta_k|\to\pi$. Since $\xi_k$ is the incenter of $\ga_k C$, there is an antipode of $\zeta_k$ in $\ga_k C$ and therefore, by Lemma~\ref{lem_regularantipodes}, we find a round $n$-sphere $\si_k$ which contains $\ga_k C$ and $\hat\ga_k\hat C$. Let $f_k:S^n\to\tits X$
be an isometric embedding with image $\si_k$. After passing to a subsequence, we may assume that $f_k$ converges pointwise with respect to the cone topology 
to a $1$-Lipschitz map $f: S^n\to\tits X$. The image $\si$ of $f$ contains $C'$ and $\hat C'$. Since $\xi$ lies in the interior of $C'$
and $\hat\xi\in\hat C'$ is antipodal to $\xi$, we see that $f$ is an isometric embedding. By Lemma~\ref{lem_regularantipodes},
the antipode of $\xi_k$ in $\si_k$ has to be the incenter of a chamber $\tilde C_k\subset\si_k$. 
Since $\zeta_k\in\hat\ga_k \hat C$ has an antipode in $\ga_k C$, we see $\hat\ga_k \hat C=\tilde C_k$. Thus the antipode of $\xi_k$
in $\si_k$ is $\hat\xi_k$, the incenter of $\hat\ga_k \hat C$, and therefore $\hat\xi_k\to\hat\xi$ as claimed.
\qed

\subsection{Induced submetries and the dimension of the base}

For this section we assume that $X$ is a non-flat locally compact CAT(0) space 
with a geometric group action $\Ga\acts X$ and where every geodesic lies in an $(n+1)$-flat.
In particular, $X$ is geodesically complete. Moreover, we assume that $X$ contains a periodic Morse $(n+1)$-flat $\hat F$ with $\hat\si:=\tits \hat F$.
By Proposition~\ref{prop_invarantset}, this means that the ideal boundary $\geo X$
contains a non-trivial proper closed symmetric $\Ga$-invariant subset. Thus, by Proposition~\ref{prop_submetry}, we have a non-trivial submetry
\[\delta:\tits X\to \Delta.\]

Our main goal is:

\bprop\label{prop_dim}
Let $X$ be a non-flat locally compact CAT(0) space where every geodesic lies in an $(n+1)$-flat.
Suppose that $X$ contains a periodic Morse $(n+1)$-flat $\hat F$ with $\hat\si:=\tits \hat F$.
Moreover, assume $\dim(\tits X)=m$ and that there exists a family of round $m$-spheres in $\tits X$ which are dense
with respect to the cone topology.
If $\tits X$ is irreducible,
then $\Delta$ is an $n$-dimensional spherical orbifold. 
\eprop

Our starting point is the observation that $\delta$  restricts to a submetry on certain round spheres in $\tits X$, including top-dimensional ones and
Tits boundaries of Morse flats.

\blem\label{lem_restrsub}
Let $\si\subset\tits X$ be a round sphere such that every direction at a point $\xi\in\si$ has an antipode in 
$\Si_\xi\si$. Then the restriction $\delta_\si:=\delta|_\si$ is a submetry $\delta_\si:\si\to\Delta$ and every $\delta_\si$-horizontal geodesic is
also $\delta$-horizontal. In particular, this applies if $\dim(\si)=\dim(\tits X)$ or if $\si$ is the Tits boundary of a periodic Morse flat.
\elem

\proof
By assumption, every point in $\tits X$ has an antipode in $\si$. Thus $\delta_\si$ is surjective.
Intersecting the fibers of $\delta$ with $\si$ clearly gives closed symmetric subsets in $\si$.
It only remains to show that they are equidistant. Let $\xi$ be a point in $\si$ and $\xi\eta$
 a $\delta$-horizontal geodesic in $\tits X$. Then we can extend $\eta\xi$ up to an antipode $\hat\eta\in\si$.
Further, denote by $\eta'$ the antipode of $\hat\eta$ in $\si$. Then $\xi\eta'$ is $\delta$-horizontal
and realizes the distance from $\xi$ to the $\delta$-fiber through $\eta$. This shows the claim.
The supplement follows from \cite[Lemma~2.1]{BL_building} and Corollary~\ref{cor_morseantipodes}.
\qed
\medskip

Let $\hat\delta:\hat\si\to\Delta$ be the induced submetry obtained from restricting $\delta$.
To get control on $\delta$, our strategy is to study $\hat\delta$ first.

We denote 
by $\{\si^1,\ldots,\si^s\}$  a set of singular hyperspheres in $\hat\si$  provided by Corollary~\ref{cor_regpoly}.
For every regular  point 
$\xi\in\hat\si_{reg}$ holds that $C_\xi$ agrees with the closure of one of the components of
$\hat\si\setminus\bigcup_{i=1}^s\si^i$. Note that $s$ will denote the number of singular hyperspheres.
The upshot is that  the hyperspheres $\si^i\subset\hat\si$ act as ``mirrors'' for horizontal geodesics.
This allows us to understand $\hat\delta$.

\blem\label{lem_invsub}
Denote by $f_i:\hat\si\to\hat\si$ the reflection at $\si^i\subset\hat\si$.
Then for all $1\leq i\leq s$, the submetry $\hat\delta$ is invariant under $f_i$, namely $\hat\delta=\hat\delta\circ f_i$.
\elem

\proof
Denote by $\tau^\pm\subset\hat\si$ the two round hemispheres determined by $\si^i$. 
Let $\tau_i\subset\tits X$ be a round $n$-hemisphere with $\tau_i\cap\hat\si=\si^i$.
By Lemma~\ref{lem_goodhemi}, the union $\hat\si\cup\tau_i$ is a convex subset of $\tits X$ which splits off $\si^i$.
Thus  $\xi$ and $f_i(\xi)$ have a common antipode in $\tau_i$. 
\qed
\medskip

This motivates the following terminology.

\bdfn
The finitely many singular hyperspheres $\si^i\subset\hat\si$, $1\leq i\leq s$,
will be called {\em reflecting hyperspheres}.
\edfn

Next we show that horizontal geodesics in $\hat\si$ can only branch off at a point on a reflecting hypersphere. 
This requires the following auxiliary lemma.

\blem\label{lem_visantiflat}
Let $\tilde\xi\in\hat\si_{reg}$ be a regular point. Further, let $\xi\in C_{\tilde\xi}$
and $\hat\xi\in\Ant_{vis}(\xi)$. Then there exists an $(n+1)$-flat $F\subset X$
which contains the set $C_{\tilde\xi}$ and the point $\hat\xi$ in its ideal boundary $\geo F$.
\elem

\proof
Let $c\subset X$ be a complete geodesic with $\geo c=\{\xi,\hat\xi\}$.
Choose a sequence  $(x_k)_{k\in\N}$ of points in $c$ with $x_k\to\hat\xi$.
Let $F_k\subset X$ be an $(n+1)$-flat which contains the geodesic ray $x_k\tilde\xi$.
By definition, $C_{\tilde\xi}\subset\geo F_k$. In particular, the geodesic ray $x_k\xi$
lies in $F_k$. By local compactness, we can extract a limit flat $F$ with the desired properties. 
\qed

\blem\label{lem_nohorbranch}
Suppose that $\xi\eta\subset\hat\si$ is a piecewise horizontal geodesic which branches at the point $\eta$.
Then there exists a reflecting hypersphere $\si^{i_0}$ which contains $\eta$ and such that $\xi\eta$ is transversal to $\si^{i_0}$, 
i.e. not contained in $\si^{i_0}$.
\elem

\proof
After possibly changing $\xi$ to an interior point of $\xi\eta$, we may assume that $\xi\eta$ is contained in $C_{\tilde\xi}$ for some $\tilde\xi\in\hat\si_{reg}$.
If the segment $\xi\eta$ branches at $\eta$, then we can extend it to a piecewise horizontal geodesic $\ga$ joining $\xi$ to an antipode $\hat\xi$
such that $\ga$ intersects $\hat\si$ precisely in $\xi\eta$. Now suppose for contradiction that the claim fails.
Since $C_{\tilde\xi}$ is a spherical polytope bounded by the hyperspheres $\si^i$, this can only happen if we can extend $\xi\eta$
inside $C_{\tilde\xi}$ to a point $\eta'$. By Lemma~\ref{lem_dense}, we find a sequence $\hat\xi_k$ of visual antipodes of $\xi$
which converge to $\hat\xi$ in the cone topology. From Lemma~\ref{lem_visantiflat}, we obtain $(n+1)$-flats $F_k$
with $C_{\tilde\xi}\subset\geo F_k$ and $\hat\xi_k\in\geo F_k$. In particular, the geodesic $\ga_k\subset\tits F_k$ from $\xi$ 
to $\hat\xi$ through the point $\eta$ also contains the point $\eta'$. Passing to a limit with respect to the cone topology, we obtain a geodesic $\ga_\infty\subset\tits X$
from $\xi$ to $\hat\xi$ which still contains the segment $\xi\eta'$.
Since $\ga$ contains $\xi\eta$, we conclude $\ga=\ga_\infty$. Contradiction. 
\qed
\medskip

An immediate consequence is the following.

\bcor\label{cor_hordirtangent}
Suppose that $\xi\in\hat\si$ lies in a single reflecting hypersphere, say $\xi\in\si^1\setminus\bigcup_{i\neq 1}\si^i$.
If a horizontal direction $v\in H_\xi$ has an antipode in $\Si_\xi\si^1$, then $v\in\Si_\xi\si^1$.
\ecor

\bcor\label{cor_nonintsplit}
Suppose that $\si:=\bigcap_{i=1}^s\si^i\neq\emptyset$. Then $\tits X$ splits non-trivially as a spherical join.
\ecor

\proof
Pick a point $\xi\in\si$. By Lemma~\ref{lem_restrsub}, $\delta$ restricts to a submetry $\hat\delta:\hat\si\to\Delta$
such that $\hat\delta$-horizontal directions are $\delta$-horizontal.
In particular, there exists a $\delta$-horizontal direction $v\in\Si_\xi\hat\si$.
By Lemma~\ref{lem_nohorbranch}, the piecewise $\delta$-horizontal geodesic $\ga$ in direction $v$
cannot branch since no $\si^i$ is transversal to $\ga$. The claim follows from Lemma~\ref{lem_trivialsplit}.
\qed
\medskip

%
%
%
%

\proof[Proof of Proposition~\ref{prop_dim}]
Since $\tits X$ is irreducible, Corollary~\ref{cor_nonintsplit}
implies $\bigcap_{i=1}^s\si^i=\emptyset$. Thus after renumbering, we may assume $\bigcap_{i=1}^{n+1}\si^i=\emptyset$.
By Lemma~\ref{lem_invsub} and Proposition~\ref{prop_splitsub}, either every geodesic in $\hat\si$ is piecewise $\hat\delta$-horizontal,
or there exists a disjoint decomposition $I^-\dot\cup I^+=\{1,\ldots,n+1\}$ and a splitting $\hat\si=\si^-\circ\si^+$
where $\si^\pm=\bigcap_{i\in I^\pm}\si^i$ such that every geodesic $\xi^-\xi^+$ is piecewise $\hat\delta$-horizontal.
By Proposition~\ref{prop_split} and Theorem~\ref{thm_split}, the latter case implies that $X$ and therefore $\tits X$ is reducible.
Thus the claim follows from Lemma~\ref{lem_polytope}.
\qed
\medskip

\subsection{Dimension of the Tits boundary and rigidity}

The aim of this section is to prove the \hyperref[thm_mainB]{Main Theorem}.
We begin with some auxiliary observations.

\bcor\label{cor_midpointsplit}
Suppose that there exists a geodesic $\ga\subset\hat\si$ of length $\pi$ whose interior intersects only one of the reflecting hyperspheres.
Moreover, assume that the intersection is perpendicular at the midpoint.
Then $\geo X$ contains a proper subset which is closed convex and symmetric.
\ecor

\proof
Let's assume that $\ga$ intersects $\si^1$ at its midpoint.
By assumption
$\D\ga\subset\bigcap_{i=2}^s\si^i$. Thus $\hat\si=\si^+\circ\si^-$ with $\si^+=\D\ga$ and $\si^-=\si^1$.
Thus the claim follows from Proposition~\ref{prop_split}.
\qed
\medskip

The following is an immediate consequence of the chain rule and the fact that the fibers of $\delta$
are symmetric.

\blem\label{lem_chainrule}
Let $\ga:S^1\to \tits X$ be an isometric embedding
and set $\xi^\pm:=\ga(\pm 1)$. Then 
\[d\delta(\xi^+)(\dot\ga(1))=d\delta(\xi^-)(\dot\ga(-1)).\]
\elem


Now we can provide
\proof[Proof of Main Theorem]
By Theorem~\ref{thm_rigidity}, it is enough to show that either $\tits X$ splits as a spherical join or $\dim(\tits X)=n$.
 
Since $\Ga$ acts geometrically, $\dim(\tits X)=m$ is finite and round $m$-spheres form a dense family in $\geo X$ (Corollary~\ref{cor_morsecyclesdense}). Let
$\bar\si\subset\tits X$ be a round $m$-sphere.
From Lemma~\ref{lem_restrsub} we obtain  induced submetries $\bar\delta:\bar\si\to\Delta$ and $\hat\delta:\hat\si\to\Delta$  by restriction. Moreover,
$\bar\delta$-horizontal directions and $\hat\delta$-horizontal directions are $\delta$-horizontal.
 From Proposition~\ref{prop_dim} we know that either $\tits X$ is reducible or $\dim(\Delta)=\dim(\hat\si)=n$ and $\Delta$ is a spherical orbifold.
Now let us assume that $\Delta$ is an $n$-dimensional spherical orbifold and $n<m$ holds.

Then, by Lemma~\ref{lem_focal}, there exists a point $\bar\xi\in\bar\si$ with $\delta(\bar\xi)=y\in\Delta^{n-1}$ 
such that the $\bar\delta$-horizontal space $\bar H_{\bar\xi}\subset\Si_{\bar\xi}\bar\si$ is a round $k$-sphere with $k\geq n$.
Let $\hat\xi\in\hat\si$ be an antipode of $\bar\xi$. In particular, $\delta(\hat\xi)=\delta(\bar\xi)=y$.

By Lemma~\ref{lem_fold}, there exists splittings $\bar H_{\bar\xi}\cong \bar L\circ\bar L^\perp$ and
$\hat H_{\hat\xi}\cong \hat K\circ\hat K^\perp$ such that the differentials $d\bar\delta(\bar\xi)$ and $d\hat\delta(\hat\xi)$ map 
$\bar L$ and $\hat K$ isometrically to $\Si_y\Delta^{n-1}\cong S^{n-2}$ and send $\bar L^\perp$ and $\hat K^\perp$ to $w$,
where $\Si_y\Delta\cong\Si_y\Delta^{n-1}\circ\{w\}$.

\begin{center}
\includegraphics[scale=0.5,trim={1cm 0cm 0cm 1cm},clip]{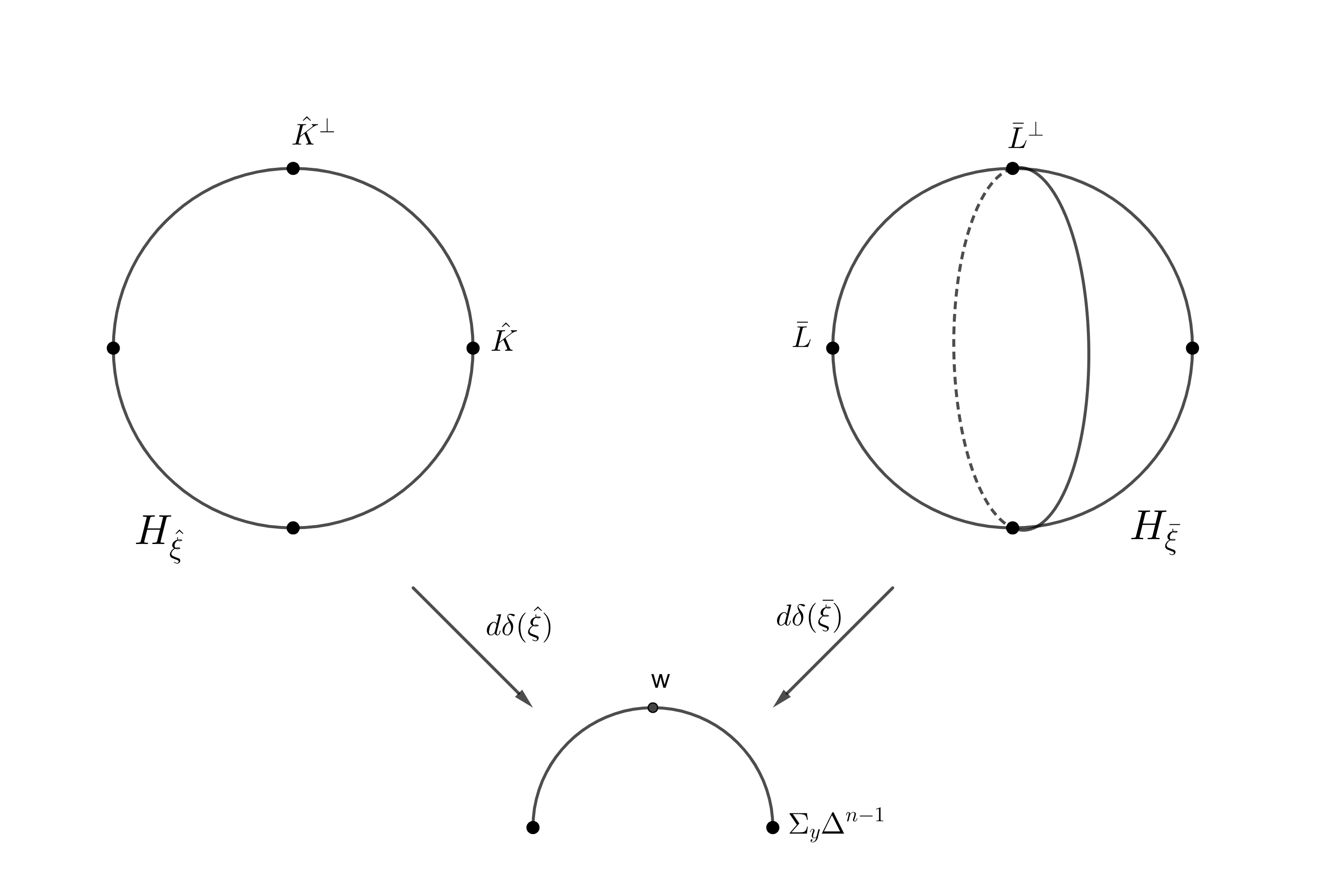}
\end{center}

Note that $\dim(\bar L^\perp)\geq 1$ while $\hat K^\perp=\{\hat v^-,\hat v^+\}$.
Let $\si\subset\tits X$ be the unique round $(k+1)$-sphere spanned by the horizontal sphere $\bar H_{\bar\xi}\subset H_{\bar\xi}$
which contains the points $\hat\xi$ and $\bar\xi$. By Lemma~\ref{lem_chainrule} the splitting of $\bar H_{\bar\xi}$ induces a splitting of the horizontal sphere
$\Si_{\hat\xi}\si\cong \bar K\circ\bar K^\perp$ such that $d\delta(\hat\xi)$ maps 
$\bar K$ isometrically to $\Si_y\Delta^{n-1}\cong S^{n-2}$ and sends $\bar K^\perp$ to $w$.

\bslem
We have $\hat K=\bar K$.
\eslem

\proof
Since $\dim(\Si_{\hat\xi}\si)=k\geq n$ and since every direction in $\Si_{\hat\xi}\si$ has an antipode in $\Si_{\hat\xi}\hat\si$, 
there is horizontal branching at $\hat\xi$.
Thus, by Lemma~\ref{lem_nohorbranch},  $\hat\xi$ has to lie on a reflecting hypersphere, say $\si^1$. 
Since $\hat\delta$ is invariant under reflections along the hyperspheres $\si^i$, $1\leq i\leq s$, it follows 
$\hat K=\Si_{\hat\xi}\si^1$. Moreover,
$\hat\xi$ cannot lie in a second reflecting hypersphere $\si^i$, $i\neq 1$, otherwise the tangent space $T_y\Delta$ 
could not split off $\R^{n-1}$, but by assumption $\delta(\hat\xi)=y\in\Delta^{n-1}$.
By Corollary~\ref{cor_morseantipodes}, every direction $\bar v\in\bar K$ has an antipode $\hat v$ in $\Si_{\hat\xi}\hat\si$.
Since $d\delta(-\hat v)=d\delta(\bar v)\in \Si_y\Delta^{n-1}$, we see $\hat v\in \hat K$.
From Corollary~\ref{cor_hordirtangent} we infer $\bar v=\hat v$.
This confirms the claim.
\qed

We conclude that every direction $\bar v'\in\bar K^\perp$ spans a round $n$-hemisphere $\bar\tau\subset\tits X$ with
$\bar\tau\cap\hat\si=\D\bar\tau$. In particular, by Lemma~\ref{lem_goodhemi}, $\angle(\bar v',\hat v^\pm)=\pi$, i.e. $\bar K^\perp\subset\Ant(\hat v^\pm)$.
Since $\dim(\bar K^\perp)\geq 1$, Lemma~\ref{lem_split} implies that a geodesic in direction $\hat v^\pm$
can only branch at times which are multiples of $\frac{\pi}{2}$. 
Since $v^\pm$ is perpendicular to $\si_1$, Corollary~\ref{cor_midpointsplit} and Theorem~\ref{thm_split}
imply that $\tits X$ is reducible. This completes the proof.
\qed

\proof[Proof of Corollary~\ref{cor_mainA}]
If $M$ is a Hadamard manifold with a geometric group action and of rank $n\geq 2$, then every geodesic lies in an $n$-flat \cite[Section~IV.4]{ballmannbook}.
Moreover, a dense set of the unit sphere bundle is tangent to periodic Morse $n$-flats \cite[Theorem~4.8]{BBS_structure}.
Thus the claim follows from the \hyperref[thm_mainB]{Main Theorem}.
\qed

%
%
%

\bibliographystyle{alpha}
\bibliography{rr_II}

\begin{thebibliography}{HKS22b}

\bibitem[AKP19]{AKP}
S.~Alexander, V.~Kapovitch, and A.~Petrunin.
\newblock Alexandrov geometry: preliminary version no. 1.
\newblock {\em arXiv:1903.08539}, 2019.

\bibitem[Bal85]{B_higher}
W.~Ballmann.
\newblock Nonpositively curved manifolds of higher rank.
\newblock {\em Ann. of Math. (2)}, 122(3):597--609, 1985.

\bibitem[Bal95]{ballmannbook}
W.~Ballmann.
\newblock {\em Lectures on spaces of nonpositive curvature}, volume~25 of {\em
  DMV Seminar}.
\newblock Birkh\"auser Verlag, Basel, 1995.
\newblock With an appendix by M. Brin.

\bibitem[Bal04]{Ballmann}
W.~Ballmann.
\newblock On the geometry of metric spaces.
\newblock {\em Preprint, lecture notes,
  http://people.mpim-bonn.mpg.de/hwbllmnn/archiv/sin40827.pdf}, 2004.

\bibitem[BB95]{BaBr_orbi}
W.~Ballmann and M.~Brin.
\newblock Orbihedra of nonpositive curvature.
\newblock {\em Inst. Hautes \'{E}tudes Sci. Publ. Math.}, (82):169--209 (1996),
  1995.

\bibitem[BB99]{BB_diam}
W.~Ballmann and M.~Brin.
\newblock Diameter rigidity of spherical polyhedra.
\newblock {\em Duke Math. J.}, 97(2):235--259, 1999.

\bibitem[BBS85]{BBS_structure}
W.~Ballmann, M.~Brin, and R.~Spatzier.
\newblock Structure of manifolds of nonpositive curvature. {II}.
\newblock {\em Ann. of Math. (2)}, 122(2):205--235, 1985.

\bibitem[Ber87]{Ber_sub}
V.~N. Berestovski\u{\i}.
\newblock ``{S}ubmetries'' of three-dimensional forms of nonnegative curvature.
\newblock {\em Sibirsk. Mat. Zh.}, 28(4):44--56, 224, 1987.

\bibitem[BGP92]{BGP}
Yu. Burago, M.~Gromov, and G.~Perelman.
\newblock A. {D}. {A}leksandrov spaces with curvatures bounded below.
\newblock {\em Uspekhi Mat. Nauk}, 47(2(284)):3--51, 222, 1992.

\bibitem[BH99]{BH}
M.~Bridson and A.~Haefliger.
\newblock {\em Metric spaces of non-positive curvature}, volume 319 of {\em
  Grundlehren der Mathematischen Wissenschaften}.
\newblock Springer-Verlag, Berlin, 1999.

\bibitem[BL06]{BL_building}
A.~Balser and A.~Lytchak.
\newblock Building-like spaces.
\newblock {\em J. Math. Kyoto Univ.}, 46(4):789--804, 2006.

\bibitem[BMS16]{BMS_affine}
H.~Bennett, C.~Mooney, and R.~Spatzier.
\newblock Affine maps between {$\rm CAT(0)$} spaces.
\newblock {\em Geom. Dedicata}, 180:1--16, 2016.

\bibitem[BS87]{BS_higher}
K.~Burns and R.~Spatzier.
\newblock Manifolds of nonpositive curvature and their buildings.
\newblock {\em Inst. Hautes \'{E}tudes Sci. Publ. Math.}, (65):35--59, 1987.

\bibitem[Ebe88]{E_sym}
P.~Eberlein.
\newblock Symmetry diffeomorphism group of a manifold of nonpositive curvature.
\newblock {\em Trans. Amer. Math. Soc.}, 309(1):355--374, 1988.

\bibitem[EH90]{EH_diff}
P.~Eberlein and J.~Heber.
\newblock A differential geometric characterization of symmetric spaces of
  higher rank.
\newblock {\em Inst. Hautes \'{E}tudes Sci. Publ. Math.}, (71):33--44, 1990.

\bibitem[GS13]{GuSw_trans}
Dan~P. Guralnik and Eric~L. Swenson.
\newblock A `transversal' for minimal invariant sets in the boundary of a
  {CAT}(0) group.
\newblock {\em Trans. Amer. Math. Soc.}, 365(6):3069--3095, 2013.

\bibitem[HKS22a]{HKS_I}
J.~Huang, B.~Kleiner, and S.~Stadler.
\newblock Morse quasiflats {I}.
\newblock {\em J. Reine Angew. Math.}, 784:53--129, 2022.

\bibitem[HKS22b]{HKS_II}
J.~Huang, B.~Kleiner, and S.~Stadler.
\newblock {M}orse {Q}uasiflats {II}.
\newblock {\em arXiv:2003.08912}, 2022.

\bibitem[HS98]{HuSch_tits}
C.~Hummel and V.~Schroeder.
\newblock Tits geometry associated with {$4$}-dimensional closed real-analytic
  manifolds of nonpositive curvature.
\newblock {\em J. Differential Geom.}, 48(3):531--555, 1998.

\bibitem[KL97]{KleinerLeeb}
B.~Kleiner and B.~Leeb.
\newblock Rigidity of quasi-isometries for symmetric spaces and {E}uclidean
  buildings.
\newblock {\em Inst. Hautes \'{E}tudes Sci. Publ. Math.}, 86:115--197, 1997.

\bibitem[KL22]{KaLy_sub}
V.~Kapovitch and A.~Lytchak.
\newblock Structure of submetries.
\newblock {\em Preprint, arXiv:2007.01325, to appear in Geom. Topol.}, 2022.

\bibitem[Kle99]{Kleiner}
B.~Kleiner.
\newblock The local structure of spaces with curvature bounded above.
\newblock {\em Math. Z.}, 231:409--456, 1999.

\bibitem[Kra11]{Kramer}
L.~Kramer.
\newblock On the local structure and the homology of {${\rm CAT}(\kappa)$}
  spaces and {E}uclidean buildings.
\newblock {\em Adv. Geom.}, 11(2):347--369, 2011.

\bibitem[Lan20]{La_orbi}
C.~Lange.
\newblock Orbifolds from a metric viewpoint.
\newblock {\em Geom. Dedicata}, 209:43--57, 2020.

\bibitem[Lee00]{Leeb}
B.~Leeb.
\newblock {\em A characterization of irreducible symmetric spaces and
  {E}uclidean buildings of higher rank by their asymptotic geometry}, volume
  326 of {\em Bonner Mathematische Schriften [Bonn Mathematical Publications]}.
\newblock Universit\"{a}t Bonn, Mathematisches Institut, Bonn, 2000.

\bibitem[LN19]{LN_gcba}
A.~Lytchak and K.~Nagano.
\newblock Geodesically complete spaces with an upper curvature bound.
\newblock {\em Geom. Funct. Anal.}, 29(1):295--342, 2019.

\bibitem[LW16]{LyWi_foliations}
A.~Lytchak and B.~Wilking.
\newblock Riemannian foliations of spheres.
\newblock {\em Geom. Topol.}, 20(3):1257--1274, 2016.

\bibitem[Lyt02]{Ly_diss}
A.~Lytchak.
\newblock {\em Allgemeine {T}heorie der {S}ubmetrien und verwandte
  mathematische {P}robleme}, volume 347 of {\em Bonner Mathematische Schriften
  [Bonn Mathematical Publications]}.
\newblock Universit\"{a}t Bonn, Mathematisches Institut, Bonn, 2002.
\newblock Dissertation, Rheinische Friedrich-Wilhelms-Universit\"{a}t Bonn,
  Bonn, 2001.

\bibitem[Lyt05]{Ly_rigidity}
A.~Lytchak.
\newblock Rigidity of spherical buildings and joins.
\newblock {\em Geom. Funct. Anal.}, 15(3):720--752, 2005.

\bibitem[Lyt10]{Ly_resol}
Alexander Lytchak.
\newblock Geometric resolution of singular {R}iemannian foliations.
\newblock {\em Geom. Dedicata}, 149:379--395, 2010.

\bibitem[Mos73]{Mos_rigidity}
G.~D. Mostow.
\newblock {\em Strong rigidity of locally symmetric spaces}.
\newblock Annals of Mathematics Studies, No. 78. Princeton University Press,
  Princeton, N.J.; University of Tokyo Press, Tokyo, 1973.

\bibitem[MR20]{MR_lap}
R.~A.~E. Mendes and M.~Radeschi.
\newblock Laplacian algebras, manifold submetries and the inverse invariant
  theory problem.
\newblock {\em Geom. Funct. Anal.}, 30(2):536--573, 2020.

\bibitem[OR22]{OR_int}
P.~Ontaneda and R.~Ricks.
\newblock {I}ntrinsic {R}ank in {CAT}(0) {S}paces.
\newblock {\em arXiv:2205.04377}, 2022.

\bibitem[OT99]{OT_cba}
Y.~Otsu and H.~Tanoue.
\newblock The riemannian structure of alexandrov spaces with curvature bounded
  above.
\newblock {\em preprint}, 1999.

\bibitem[PR72]{PR_cartan}
G.~Prasad and M.~S. Raghunathan.
\newblock Cartan subgroups and lattices in semi-simple groups.
\newblock {\em Ann. of Math. (2)}, 96:296--317, 1972.

\bibitem[PS09]{PaSw_bound}
P.~Papasoglu and E.~Swenson.
\newblock Boundaries and {JSJ} decompositions of {CAT}(0)-groups.
\newblock {\em Geom. Funct. Anal.}, 19(2):559--590, 2009.

\bibitem[Sta22a]{St_rr}
S~Stadler.
\newblock {CAT}(0) spaces of higher rank {I}.
\newblock {\em Preprint}, 2022.

\bibitem[Sta22b]{St_rrIII}
S.~Stadler.
\newblock {R}ank {Rigidity} for {CAT}(0) spaces without {3}-flats.
\newblock {\em Preprint}, 2022.

\bibitem[Tho22]{Thor}
Gudlaugur Thorbergsson.
\newblock From isoparametric submanifolds to polar foliations.
\newblock {\em S\~{a}o Paulo J. Math. Sci.}, 16(1):459--472, 2022.

\end{thebibliography}

\Addresses

\end{document}